\documentclass[12pt]{amsart}
\newtheorem{thm}{Theorem}[section]

\newtheorem{lem}[thm]{Lemma}
\newtheorem{prop}[thm]{Proposition}
\theoremstyle{definition}

\theoremstyle{remark}

\numberwithin{equation}{section}

\newcommand{\Real}{\mathbb R}
\newcommand{\Sphere}{\mathbb S}
\newcommand{\eps}{\varepsilon}

\begin{document}
\vskip3truecm

\centerline{\bf{ON THE SPECTRUM OF THE
LAPLACE-BELTRAMI}}\centerline{\bf{ OPERATOR
 FOR $p$-FORMS FOR A CLASS OF}}\centerline{\bf{ WARPED PRODUCT
METRICS }} \par

\vskip3truecm \centerline{Francesca Antoci}

\vskip3truecm

\noindent Affiliation: Dipartimento di Matematica, Politecnico di
Torino, C.so Duca degli Abruzzi 24, 10129,  Torino, Italy. \par
\vskip1truecm \noindent e-mail address: antoci@calvino.polito.it
\vskip1truecm \noindent Address for correspondence: Dipartimento
di Scienze Matematiche, Universit\`a degli Studi di Trieste, Via
Valerio 12/b, 34127, Trieste, Italy.

\thanks{}%

\begin{abstract} We explicitely compute the essential spectrum of the \\
Laplace-Beltrami operator for $p$-forms for the class of warped product
metrics $d\sigma^2= y^{2a}dy^2 + y^{2b}d\theta_{\partial M}^2$, where $y$
is a boundary defining function on a compact manifold with boundary $M$.
\end{abstract}

\newpage
\section{Introduction}
The spectrum of the Laplace-Beltrami operator on complete noncompact
Riemannian manifolds in its connections with the geometrical properties
 has been investigated by many authors; however, they
have mainly studied the case of scalar functions, whilst less is
known about $p$-forms, since the complicated local expression of
the Laplace-Beltrami operator on $p$-forms makes difficult any
explicit computation. Hence, the attention has mainly focused on
particular classes of Riemannian manifolds in which these
difficulties can be bypassed thanks to the presence of symmetries
or ``asymptotic symmetries". This is the case for the hyperbolic
space $\mathbb H^n$ (\cite{Donnelly}), or, more generally, for
rotationally symmetric Riemannian metrics (\cite{Eichhorn}) or
manifolds with cylinderlike ends
(\cite{Eichhorn2},\cite{Eichhorn3}), where a decomposition
technique introduced by Dodziuk (\cite{Dodziuk}) considerably
simplifies the problem. By this technique Donnelly
(\cite{Donnelly}) computed the spectrum of the Laplace-Beltrami
operator on the hyperbolic space $\mathbb H^n$, and Eichhorn
(\cite{Eichhorn}, \cite{Eichhorn2}, \cite{Eichhorn3}) obtained
results on the spectrum of the Laplace-Beltrami operator for
$p$-forms on  Riemannian manifolds with cylindrical ends $M=
(0,+\infty)\times N$ endowed with the Riemannian metric
$d\sigma^2= dt^2+g(t)d\theta^2$, where $t\in (0,+\infty)$,
$g(t)>0$ and $N$ is a compact manifold endowed with the Riemannian
metric $d\theta^2$. \par In the present paper, we compute the
essential spectrum of the Laplace-Beltrami operator for $p$-forms
for a class of warped product metrics. Namely, let $\overline M$
be a compact Riemannian $n$-dimensional manifold with boundary,
and let $y$ be a boundary-defining function; we endow the interior
$M$ of $\overline M$ with a Riemannian metric $ds^2$ such that in
a small tubular neighbourhood of $\partial M$ in $M$, $ds^2$ takes
the form
\begin{equation}\label{form0} ds^2=
e^{-2(a+1)t}\,dt^2+e^{-2bt}d\theta^2_{\partial M},\end{equation}
where $t:=-\log y\in (c,+\infty)$ and $d\theta_{\partial M}^2$ is the
Riemannian metric on $\partial M$. For $a\leq -1$ (see \cite{Melrose}) the
Riemannian metric (\ref{form0}) is complete. Hence, the Laplace-Beltrami
operator is essentially selfadjoint on the smooth, compactly supported
$p$-forms; we compute the essential spectrum of its closure $\Delta^p_M$
in dependence on the parameters $a$ and $b$.\par
Since Eichhorn (\cite{Eichhorn}) showed that, for any $p\in [0,n]$,
the essential
spectrum of
$\Delta^p_M$ coincides with the essential spectrum of the
Friedrichs extension $(\Delta_M^p)^F$ of $\Delta^p_M$ on smooth $p$-forms
with compact support in $(c,+\infty)\times \partial M$, we compute the
essential
spectrum of $(\Delta_M^p)^F$. To this purpose we use an orthogonal
decomposition of $L^2_p((c,+\infty)\times \partial M)$ analogue to those
employed by Eichhorn and Donnelly. The decomposition is carried out in two
steps: first,
thanks to the Hodge decomposition on $\partial M$, we write any $p$-form
$\omega$ on $(c,+\infty)\times \partial M$ as
\begin{equation}\label{ddec0}\omega= \omega_{1\delta}\oplus
\omega_{2d}\wedge dt \oplus (\omega_{1d}\oplus
\omega_{2\delta}\wedge dt),\end{equation} where $\omega_{1\delta}$
(resp. $\omega_{1d}$) is a coclosed (resp. closed) $p$-form on
$\partial M$ parametrized by $t$ and $\omega_{2\delta}$ (resp.
$\omega_{2d}$) is a coclosed (resp. closed) $(p-1)$-form on
$\partial M$ parametrized by $t$. This yields the orthogonal
decomposition \begin{multline}\label{decL0} L^2_p((c,+\infty)\times
\partial M)= \\ {\mathcal L
}_1((c,+\infty)\times \partial M)\oplus {\mathcal L}_2
((c,+\infty)\times \partial M)\oplus {\mathcal L}_3((c,+\infty)\times
\partial M),\end{multline}
where for any $\omega\in L^2_p((c,+\infty)\times \partial M)$
$$\omega_{1\delta}\in {\mathcal L}_1((c,+\infty)\times \partial
M),$$ $$\omega_{2d}\wedge dt\in {\mathcal L}_2((c,+\infty)\times
\partial M),$$ and $$(\omega_{1d}\oplus \omega_{2\delta }\wedge
dt)\in {\mathcal L}_3((c,+\infty)\times \partial M).$$ Since
$\Delta^p_M$ is invariant under (\ref{decL0}), we get the
corresponding decomposition of $(\Delta_M^p)^F$ $$
(\Delta_M^p)^F=(\Delta_{M1}^p)^F \oplus (\Delta_{M2}^p)^F \oplus
(\Delta_{M3}^p)^F,$$ whence $$ \sigma_{\rm
ess}(\Delta^p_M)=\sigma_{\rm ess}((\Delta_M^p)^F)= \bigcup_{i=1}^3
\sigma_{\rm ess}((\Delta^p_{Mi})^F).$$ Hence, the calculus of the
essential spectrum of $\Delta_M^p$ can be reduced to the
determination of the essential spectra of $(\Delta^p_{Mi})^F$, for
$i=1,2,3$. Moreover, since the Hodge $*$ operator maps
isometrically $p$-forms of ${\mathcal L}_1((c,+\infty)\times
\partial M)$ onto $(n-p)$-forms of ${\mathcal
L}_2((c,+\infty)\times \partial M)$, it suffices to consider the
cases $i=1,3$.\par
\par The second step consists in the decomposition of
$\omega_{1\delta}$ (resp. of $\omega_{2d}$, $\omega_{2\delta}$)
according to an orthonormal basis of coclosed $p$-eigenforms
(resp. closed and coclosed $(p-1)$-eigenforms) of
$\Delta^p_{\partial M}$ (resp. $\Delta_{\partial M}^{p-1}$) on $\partial
M$. In this way, up to a
unitary equivalence, the spectral analysis of $(\Delta_{Mi}^p)^F$,
$i=1,3$, can be reduced to the investigation of the spectra of a
countable number of Sturm-Liouville operators $(D_{i\lambda})^F$ on
$(c,+\infty)$, parametrized by the eigenvalues $\lambda_k^p$ of
$\Delta^p_{\partial M}$ if $i=1$, and
by the eigenvalues $\lambda^{p-1}_k$ of $\Delta^{p-1}_{\partial M}$
 if $i=3$. We remark that, since we deal with an
infinite family of operators, we only get the inclusion $$
\bigcup_{k=1}^{\infty}\sigma_{\rm
ess}((D_{i\lambda_k})^F)\subseteq \sigma_{\rm
ess}((\Delta^p_{Mi})^F).$$ In order to fully determine the
essential spectrum, we need additional estimates on the behaviour
of the isolated eigenvalues of finite multiplicity of
$(D_{i\lambda_k})^F$ as $k\rightarrow +\infty$.\par For $i=1$ the
situation is simpler; the essential spectrum of $(D_{1\lambda})^F$
can be easily computed, and when necessary we can easily deduce
the estimates on the isolated eigenvalues of finite multiplicity.
Moreover, as already noted by Eichhorn (\cite{Eichhorn2}), this
part of the spectrum depends strongly on the geometry of $M$ (more
precisely, on the asymptotic sectional curvature). \par On the
other hand, the analysis of $(\Delta_{M3}^p)^F$ is difficult,
since the operators $(D_{3\lambda})^F$ are coupled systems of
Sturm-Liouville operators on $L^2(c,+\infty)$ $\oplus
L^2(c,+\infty)$. However, the essential spectrum of
$D_{3\lambda}^F$ can be determined through perturbation theory,
taking, when necessary, as ``unperturbed operator" an operator
$D_{30}^F$ which is invariant under the transformation $(w_1\oplus
w_2)\longrightarrow (w_2\oplus w_1)$. Indeed, due to this
invariance, the essential spectrum of $D_{30}^F$ can be computed
as the union of the essential spectra of its restrictions to the
subspaces ${\mathcal V}_1=\left\{ w_1\oplus w_2 \,\mid\,w_1=w_2
\right\}$ and ${\mathcal V}_2= \left\{ w_1 \oplus w_2\,|\,w_1=-w_2
\right\}$ of $L^2(c,+\infty)\oplus L^2(c,+\infty)$.\par More
difficult is to get the estimates on the behaviour of the isolated
eigenvalues of finite multiplicity of $(D_{3\lambda_k})^F$ as
$k\rightarrow +\infty$. However, this difficulty can be overcome
observing that, roughly speaking, the differential $d^p_M \omega$
of a smooth $p$-eigenform $\omega$ of $\Delta^p_M$ in ${\mathcal
L}_3((c,+\infty)\times \partial M)$ is a smooth $(p+1)$-eigenform
of $\Delta^{p+1}_M$ in ${\mathcal L}_2((c,+\infty)\times
\partial M)$, whilst its
codifferential $\delta^p_M\omega$ is a smooth $(p-1)$-eigenform of
$\Delta^{p-1}_M$ in ${\mathcal L}_1((c,+\infty)\times \partial
M)$. This entails the existence of a link between the positive
parts of the spectra for forms of different types and degrees.
Namely, one can prove that if there exists a sequence
$\left\{\mu_k\right\}$ of positive eigenvalues of
$(\Delta_{M3}^p)^F$ and a corresponding sequence of eigenforms
$\left\{\Phi_k\right\}$ such that, for every $k\in \mathbb N$,
$(\Delta_{M3}^p)^F\Phi_k-\mu_k \Phi_k=0$ and $\mu_k \rightarrow
\mu\geq 0$ as $k\rightarrow +\infty$, then either $\mu$ is in the
essential spectrum of $(\Delta_{M1}^{(p-1)})^F$ or $\mu$ is in the
essential spectrum of $(\Delta_{M2}^{(p+1)})^F$. This observation,
combined with the information which we already have on the
essential spectra of $(\Delta_{M1}^p)^F$ and $(\Delta_{M2}^p)^F$
for every degree $p$, allows to compute completely $\sigma_{\rm
ess}(\Delta^p_M)\setminus \left\{0\right\}$.\par However, in this
way it is not possible to determine if $0$ lies or not in the
essential spectrum; indeed, we cannot decide whether $0$ is an
eigenvalue of infinite multiplicity of $(\Delta_{M3}^p)^F$. This
can be established, instead, if there is a rotational symmetry,
that is, if $\overline{M}$ is the unitary ball $\overline
{B(\underline 0,1)}$ in $\Real^n$ and the Riemannian metric is
globally invariant under rotations. In fact, if the Riemannian
metric coincides with the Euclidean metric in a small
neighbourhood $B(\underline 0,\epsilon)$ of $\underline 0$, a
slight modification of a classical result of Dodziuk
(\cite{Dodziuk}) allows to check whether $0\in
\sigma_p(\Delta^p_M)$ and, if so, to determine its multiplicity.
Adding this new information, in this case we can describe
completely the essential spectrum of $\Delta^p_M$. \par Now we
discuss briefly our results. If $a=-1$, the infimum of
$\sigma_{\rm ess}(\Delta^p_M)\setminus \left\{0\right\}$ for a
given dimension $n$ depends on the degree $p$ of the forms. Since
there is a strong dependence of the spectral properties of the
operators on the sign of $b$, we consider separately the cases
$b<0$, $b=0$, $b>0$.\par If $a=-1$ and $b<0$, the situation is
similar to the hyperbolic case (see \cite{Donnelly},
\cite{Antoci}). With no rotational symmetry assumptions, $$
\sigma_{\rm ess}(\Delta^p_M)\setminus \left\{0\right\}=$$$$=
\left[ \min
\left\{\left(\frac{n-2p-1}{2}\right)^2b^2,\left(\frac{n-2p+1}{2}
\right)^2b^2 \right\}, +\infty \right),$$ whilst in the
rotationally symmetric case $$\sigma_{\rm ess}(\Delta^p_M)=
\left[\min\left\{\left(\frac{n-2p-1}{2}\right)^2 b^2,
\left(\frac{n-2p+1}{2}\right)^2b^2\right\},+\infty\right) $$ for
$p\not= \frac{n}{2}$ and $$\sigma_{\rm ess}(\Delta^p_M)=
\left\{0\right\}\cup \left[\frac{b^2}{4},+\infty \right) $$ for
$p=\frac{n}{2}$.\par If $a=-1$ and $b=0$, we get the well-known
result for cylindrical ends $$ \sigma_{\rm ess}(\Delta^p_M)=
\bigcup_{k}\left([\lambda_k^p,+\infty) \cup
[\lambda^{p-1}_k,+\infty)\right),$$ where $\lambda^p_k$ (resp.
$\lambda^{p-1}_k$) denote as usual the eigenvalues of the
Laplace-Beltrami operator $\Delta^p_{\partial M}$ (resp.
$\Delta_{\partial M}^{p-1}$) on $\partial M$.\par If $a=-1$ and
$b>0$, the essential spectrum of $\Delta^p_M$ somewhat reflects
the cohomology of $\partial M$. Indeed, if both the $p$-th and the
$(p-1)$-th Betti numbers of $\partial M$ vanish, $\sigma_{\rm
ess}(\Delta^p_M)\setminus \left\{0\right\} = \emptyset;$ if the
$p$-th Betti number of $\partial M$ is different from zero while
the $(p-1)$-th Betti number of $\partial M$ vanishes, $$
\sigma_{\rm ess}(\Delta^p_M)\setminus \left\{0\right\}= \left[
\left(\frac{n-2p-1}{2}\right)^2b^2,+\infty\right);$$
 if the $p$-th Betti number of $\partial M$ vanishes while the
$(p-1)$-th Betti number of $\partial M$ is different from zero, $$
\sigma_{\rm ess}(\Delta^p_M)\setminus \left\{ 0\right\} =
\left[\left( \frac{n-2p+1}{2}\right)^2 b^2, +\infty \right);$$
finally, if both the $p$-th and the $(p-1)$-th Betti numbers of
$\partial M$ are different from zero,  $$\sigma_{\rm
ess}(\Delta^p_M)\setminus \left\{ 0 \right\} = $$
$$=\left[\min\left\{\left(\frac{n-2p-1}{2}\right)^2b^2,\left(\frac{n-2p+1}
{2}\right)^2b^2\right\}, +\infty\right). $$ In particular, in the
rotationally symmetric case, where $\partial M=\Sphere^{n-1}$, we
find that $0\notin \sigma_{\rm ess}(\Delta^p_M)$; hence, for
$1<p<n-1$ $\sigma_{\rm ess}(\Delta^p_M)= \emptyset$, whilst for
$p\in \left\{ 0,1,n-1,n\right\}$ $\sigma_{\rm ess}(\Delta^p_M)=
\left[\left(\frac{n-1}{2}\right)^2b^2,+\infty \right).$
\par
If $a<-1$, for $b\not=0$ the infimum of $\sigma_{\rm
ess}(\Delta^p_M)\setminus \left\{0\right\}$ is $0$. Again, much
depends on the sign of $b$. \par If $b<0$, for every $p\in[0, n]$,
$\sigma_{\rm ess}(\Delta^p_M)=[0,+\infty)$.\par If $b=0$,
as before, $$\sigma_{\rm ess}(\Delta^p_M)= \bigcup_k
\left([\lambda^p_k, +\infty) \cup [\lambda^{p-1}_k,+\infty)\right)
.$$ If finally $b>0$, also for $a<-1$ the essential spectrum of
$\Delta^p_M$ reflects somehow the cohomology of the boundary $\partial M$.
As a matter of fact, if both the $p$-th and the $(p-1)$-th Betti
numbers of $\partial M$ vanish, $\sigma_{\rm
ess}(\Delta^p_M)\setminus \left\{0\right\}= \emptyset$. If, on the
contrary, at least one among the $p$-th and the $(p-1)$-th Betti
numbers of $\partial M$ is different from zero, then $\sigma_{\rm
ess}(\Delta^p_M)=[0,+\infty)$. As a consequence, in the rotationally
symmetric case, where $\partial M=\Sphere^{n-1}$, if $1<p<n-1$
$\sigma_{\rm ess}(\Delta^p_M)=\emptyset$, whilst if $p\in
\left\{0,1,n-1,n \right\}$ $\sigma_{\rm
ess}(\Delta^p_M)=[0,+\infty)$.\par
A more detailed analysis might yield a description of the
absolutely continuous spectrum of $\Delta^p_M$. However, at least for the
time being, this turns out to be very difficult because the Fourier
transform, which is a basic tool in scattering
theory, is not available in this context.\par The paper is organized as
follows: in section 2 we
introduce some preliminary facts and basic notations. In section 3
we describe in some details the decomposition techniques. The actual
calculus of the essential spectrum (with $0$ excluded) in the
general case is developed in section 4 for $a=-1$ and in section 5
for $a<-1$. In section 6, finally, we fully determine the
essential spectrum of $\Delta^p_M$ in the rotationally symmetric
case.
\section{Preliminary facts and notations}
Let $\overline M$ be a compact, $n$-dimensional manifold with boundary,
and let $y$ be a positive defining function for $N:=\partial M$:
$$ y\geq 0,\quad \quad y^{-1}(0)= N, \quad \quad dy_{|N}=0.$$
We endow the interior $M$ of $\overline M$ with a Riemannian metric which,
in a tubular neighbourhood $(0,\eps)\times N$ of $N$ in
$M$, is given by
\begin{equation}\label{Mel}d\sigma^2= y^{2a}\,dy^2 +
y^{2b}\,d\theta_N^2,\end{equation}
where $d\theta_N^2$ is a Riemannian metric on $N$.\par
Under the change of variables $t=-\log y$, (\ref{Mel}) is trasformed into
the
Riemannian metric
\begin{equation}\label{exp} ds^2= e^{-2(a+1)t}\,dt^2 +
e^{-2bt}\,d\theta^2_N,\end{equation} defined on $(c,+\infty)
\times N$, where $c=-\log \eps$.\par It is well-known (see
\cite{Melrose}) that
a Riemannian metric of this kind is complete if and only if $a\leq
-1$. Hence, throughout the paper, we will suppose that $a\leq -1$.
\par \noindent For $p=0,...,n$, we will denote by $C^{\infty}(\Lambda^p(M))$
(resp. $C^{\infty}(\Lambda^p$ $((c, +\infty)\times N))$) the space
of all the smooth $p$-forms on $M$ (resp. on $(c,+\infty)\times
N$), and by $C^{\infty}_c(\Lambda^p(M))$ (resp.
$C^{\infty}_c(\Lambda^p((c,+\infty)\times N))$) the set of all the
smooth $p$-forms with compact support in $M$ (resp. in
$(c,+\infty)\times N$).\par For any $\omega \in
C^{\infty}(\Lambda^p(M))$, we will denote by $|\omega(x)|$ the
norm induced by the Riemannian metric on the fiber over $x$, given
in local coordinates by
$$|\omega(x)|^2=g^{i_1j_1}(x)...g^{i_pj_p}(x)
\omega_{i_1...i_p}(x)\omega_{j_1...j_p}(x) ,$$ where $g^{ij}$ is
the expression of the Riemannian metric in local coordinates.
Following \cite{deRham}, we will denote by $d^p_M$, $*_M$,
$\delta^p_M$ respectively the differential, the Hodge $*$ operator
and the codifferential for $p$-forms on $M$. $\Delta^p_M$ will stand for
the
Laplace-Beltrami operator acting on $p$-forms $$ \Delta^p_M = d^{p-1}_M
\delta^p_M + \delta^{p+1}_M d^p_M,$$ which is expressed in local
coordinates
by the Weitzenb\"ock formula $$(\Delta^p_M \omega)_{i_1...i_p}= -
g^{ij}\nabla_i \nabla_j \omega_{i_1...i_p}+ \sum_j R^{\alpha}_j
\omega_{i_1...\alpha...i_p}+ \sum_{j,l\not=
j}R^{\alpha\,\beta\,}_{\,i_j\,i_l}\omega_{\alpha
i_1...\beta...i_p} ,$$ where $\nabla_i \omega$ is the covariant
derivative of $\omega$ with respect to the Riemannian metric, and
$R^i_j$, $R^{i\,j\,}_{\,k\,l}$ denote respectively the local
components of the Ricci tensor and of the Riemann tensor induced
by the Riemannian metric. \par As usual, $L^2_p(M)$ will denote
the completion of $C^{\infty}_c(\Lambda^p(M))$ with respect to the
norm $\|\omega\|_{L^2_p(M)}$ induced by the scalar product
$$\langle\omega,\tilde\omega\rangle_{L^2_p(M)}:= \int_M \omega
\wedge
*_M\tilde \omega
; $$
$\|\omega\|_{L^2_p(M)}$ can also be written as $$
\|\omega\|^2_{L^2_p(M)}= \int_M |\omega(t,\theta)|^2 dV_M,$$ where
$dV_M$ is the volume element of $(M,ds^2)$.\par It is well-known
that, since the Riemannian metric on $M$ is complete, the
Laplace-Beltrami operator is essentially selfadjoint on
$C^{\infty}_c(\Lambda^p(M))$, for $p=0,...,n$. We will denote by
the same symbol $\Delta^p_M$ its closure, and we will study its
essential spectrum.\par We recall that, given a selfadjoint
operator $A:{\mathcal D}(A)\rightarrow {\mathcal H}$ in a Hilbert
space ${\mathcal H}$, its essential spectrum is defined as the set
$$\sigma_{\rm ess}(A):=\sigma(A)\setminus \sigma_d(A), $$ where
$\sigma_d(A)$, called the discrete spectrum of $A$, is the set of
isolated eigenvalues of finite multiplicity of $A$. The essential
spectrum of $A$ can be also characterized in terms of sequences as
follows: $\mu \in \sigma_{\rm ess}(A)$ if and only if there exists
a Weyl sequence for $\mu$, that is, a sequence
$\left\{w_n\right\}\subset {\mathcal D}(A)$ with no convergent
subsequences in ${\mathcal H}$, bounded in ${\mathcal H}$ and such
that $$\lim_{n\rightarrow +\infty}(A-\mu) w_n =0 \quad\hbox{in
${\mathcal H}$}.$$ \par If $A$ is a selfadjoint operator on a
Hilbert space ${\mathcal H}$, an operator $C$ such that ${\mathcal
D}(A)\subset {\mathcal D}(C)$ is called relatively compact with
respect to $A$ if $C$ is compact from ${\mathcal D}(A)$ with the
graph norm $\|.\|_A$  given by $$\|\phi\|^2_A = \|\phi\|_{\mathcal
H}^2+ \|A\phi\|_{\mathcal H}^2$$ to ${\mathcal H}$ with the norm
$\|.\|_{\mathcal H}$. Moreover (see \cite{Reed-Simon}), given a
selfadjoint operator $A$ on a Hilbert space ${\mathcal H}$, and a
symmetric operator $C$  relatively compact with respect to $A$, if
$B=A+C$ is selfadjoint on ${\mathcal D}(A)$, then $$ \sigma_{\rm
ess}(A)= \sigma_{\rm ess}(B) .$$ For any symmetric, positive
operator $A:{\mathcal D}(A)\rightarrow {\mathcal H}$ on a Hilbert
space ${\mathcal H}$, we will denote by $A^F$ its Friedrichs
extension, defined by $${\mathcal D}(A^F)=X_A \cap {\mathcal
D}(A^*), $$ $$A^F w = A^* w ,$$ where $X_A$ denotes the completion
of ${\mathcal D}(A)$ with respect to the graph norm. \par In
particular, we will denote by $(\Delta^p_{M,c})^F$ the Friedrichs
extension of the restriction of $\Delta^p_M$ to
$C^{\infty}_c(\Lambda^p((c,+\infty)\times N))$. \par In
\cite{Eichhorn}, J. Eichhorn showed that the essential spectrum of
the Laplace-Beltrami operator acting on $p$-forms on a complete
noncompact Riemannian manifold is equal to the essential spectrum
of the Friedrichs extension of the restriction of the same
operator to the smooth $p$-forms with compact support in $M
\setminus K$, where $K\subset M$ is any compact subset of $M$. In
the present case, this yields $$\sigma_{\rm{ess}}(\Delta^p_M)=
\sigma_{\rm{ess}}((\Delta^p_{M,c})^F) ,$$ where $c >-\log \eps$
can be taken arbitrarily large. \par

\section{Hodge decomposition}
In the present section let us suppose, more generally, that the Riemannian
metric $ds^2$ in $(c,+\infty)\times N$, takes the form
\begin{equation}\label{fg}ds^2= f(t)\,dt^2 +
g(t)\,d\theta_N^2,\end{equation}
where $f(t)>0$ and $g(t)>0$ for any $t\in (c,+\infty)$. Moreover, we will
suppose that $c>0$ is fixed and we will write $(\Delta^p_M)^F$ for
$(\Delta^p_{M,c})^F$. \par
Now, any $\omega \in C^{\infty}_c(\Lambda^p((c,+\infty)\times N))$ can be
decomposed as
\begin{equation}\label{decom1} \omega= \omega_1 + \omega_2
\wedge dt,\end{equation}
where $\omega_1$ and $\omega_2$ are respectively a $p$-form and a
$(p-1)$-form on $N$ depending on $t$. An easy
computation shows that $*_M \omega$ can be expressed in terms of
decomposition
(\ref{decom1}) as
\begin{equation}\label{Hodge1}*_M \omega=
(-1)^{n-p}g^{\frac{n-2p+1}{2}}(t) f^{-\frac{1}{2}}(t)
*_N \omega_2 +
g^{\frac{n-2p-1}{2}}(t)f^{\frac{1}{2}}(t)*_N \omega_1\wedge
dt, \end{equation}
where $*_N$ denotes the Hodge $*$ operator on
$N$. Moreover, $d^p_M$ and $\delta^p_M$ split respectively as
\begin{equation}\label{dM1}d^p_M\omega= d^p_N\omega_1 + \left\{
(-1)^p \frac{\partial \omega_1}{\partial t}+ d^{p-1}_N \omega_2
\right\}\wedge dt ,\end{equation}
\begin{multline}\label{deltaM1}\delta^p_M \omega= g^{-1}(t)
\delta^p_N \omega_1 + (-1)^p f^{-\frac{1}{2}}
g^{\frac{-n-1+2p}{2}}\frac{\partial}{\partial t}
\left(f^{-\frac{1}{2}} g^{\frac{n+1-2p}{2}}\omega_2\right)+\\+
g^{-1} \delta^{p-1}_N \omega_2 \wedge dt,\end{multline} where
$d^p_N$ is the differential on $N$ and $\delta^p_N$ is the
codifferential on $N$. \par Moreover, the $L^2$-norm of $\omega
\in C^{\infty}(\Lambda^p((c,+\infty)\times N))\cap
L^2_p((c,+\infty)\times N)$ can be written as
\begin{multline}\label{norm1} \|\omega\|^2_{L^2_p((c,+\infty)\times N)}= \int_c^{+\infty}
g^{\frac{n-2p-1}{2}}(s) f^{\frac{1}{2}}(s)
\|\omega_1(s)\|^2_{L^2_p(N)} \,ds +\\
+ \int_c^{+\infty} g^{\frac{n+1-2p}{2}}(s)f^{-\frac{1}{2}}(s)
\|\omega_2(s)\|^2_{L^2_{p-1}(N)} \,ds,
\end{multline}
where $\|.\|_{L^2_p(N)}$ is the $L^2$-norm for $p$-forms on
$N$.   \par
From (\ref{dM1}) and (\ref{deltaM1}), a
lengthy but straightforward computation gives
$$ \Delta^p_M \omega= (\Delta^p_M \omega)_1 + (\Delta^p_M \omega)_2\wedge
dt,$$
with
\begin{multline}\label{LB11}(\Delta^p_M \omega)_1=
g^{-1}(t)\Delta^p_N \omega_1 + (-1)^p f^{-1}(t)
g^{-1}(t) \frac{\partial g}{\partial t} d_N^{p-1}\omega_2 +\\
-f^{-\frac{1}{2}}(t) g^{\frac{-n+1+2p}{2}}(t)
\frac{\partial}{\partial t} \left( f^{-\frac{1}{2}}(t)
g^{\frac{n-1-2p}{2}}(t) \frac{\partial \omega_1}{\partial
t} \right)
\end{multline} and
\begin{multline}\label{LB12}(\Delta_M^p \omega)_2 =
g^{-1}(t) \Delta^{p-1}_N \omega_2  + (-1)^{p}g^{-2}(t)
\frac{\partial g}{\partial t} \delta^p_N \omega_1  +\\ -
\frac{\partial}{\partial t} \left\{f^{-\frac{1}{2}}(t)
g^{\frac{-n-1+2p}{2}}(t) \frac{\partial}{\partial
t}\left(f^{-\frac{1}{2}}(t) g^{\frac{n+1-2p}{2}}(t)
\omega_2\right) \right\}, \end{multline} where $\Delta^p_N$ (resp.
$\Delta_N^{p-1}$) stands for the Laplace-Beltrami operator for
$p$-forms (resp. for $(p-1)$-forms) on $N$.\par Since for every
smooth $\omega\in L^2_p((c,+\infty)\times N)$ we have that
$\omega_1 \in L^2_p((c,+\infty)\times N)$, $\omega_2 \wedge dt \in
L^2_p((c,+\infty)\times N)$ and $$ \langle\omega_1,\omega_2\wedge
dt\rangle_{L^2_p((c,+\infty)\times N)}=0,$$ (\ref{decom1}) gives
rise to an orthogonal decomposition of $L^2_p((c,+\infty)\times
N)$ into two closed subspaces. However, (\ref{LB11}) and
(\ref{LB12}) show that $\Delta^p_M$ is not invariant under this
decomposition. As a consequence, further decompositions are
required.\par It is well-known that, for $0\leq p \leq n-1$, $$
C^{\infty}(\Lambda^p(N))= dC^{\infty}(\Lambda^{p-1}(N))\oplus
\delta C^{\infty}(\Lambda^{p+1}(N)) \oplus {\mathcal H}^p(N),$$
where ${\mathcal H}^p(N)$ is the space of harmonic $p$-forms on
$N$, and the decomposition is orthogonal in $L^2_p(N)$. Hence, for
$0\leq p \leq n-1$,
$$L^2_p(N)=\overline{dC^{\infty}(\Lambda^{p-1}(N))}\oplus
\overline{\delta C^{\infty}(\Lambda^{p+1}(N))} \oplus {\mathcal
H}^p(N). $$ As a consequence, for $1\leq p \leq n-1$, every
$\omega\in L^2_p((c,+\infty)$ $\times N)$ can be written as
$$\omega= \omega_{1\delta}\oplus \omega_{2d}\wedge dt \oplus
(\omega_{1d} \oplus \omega_{2\delta}\wedge dt),$$ where
$\omega_{1\delta}$ (resp. $\omega_{1d}$) is a coclosed (resp.
closed) $p$-form on $N$ parame-\\ trized by $t$, and
$\omega_{2\delta}$ (resp. $\omega_{2d}$) is a coclosed (resp.
closed) $(p-1)$-form on $N$ parametrized by $t$. In this way we
get the orthogonal decomposition $$ L^2_p((c,+\infty)\times N)=$$
$$= {\mathcal L}_1((c,+\infty)\times N)\oplus {\mathcal
L}_2((c,+\infty)\times N)\oplus{\mathcal L}_3((c,+\infty)\times
N),$$ where, for every $\omega \in L^2_p((c,+\infty)$ $\times N)$,
$$\omega_{1\delta}\in {\mathcal L}_1((c,+\infty)\times N),$$ $$
\omega_{2d}\wedge dt \in {\mathcal L}_2((c,+\infty)\times N)$$ and
$$ (\omega_{1d} \oplus \omega_{2\delta}\wedge dt)\in {\mathcal
L}_3((c,+\infty)\times N).$$ Since $$d^p_N \Delta^p_N
=\Delta^{p+1}_N d^p_N,\quad \quad \delta^p_N \Delta^p_N
=\Delta^{p-1}_N \delta^p_N,$$ $$\frac{\partial}{\partial
t}d^p_N=d^p_N \frac{\partial}{\partial t}, \quad \quad
\frac{\partial}{\partial t}\delta^p_N=\delta^p_N
\frac{\partial}{\partial t},$$ the Laplace-Beltrami operator is
invariant under this decomposition, and
$$(\Delta^p_M)^F=(\Delta^p_{M1})^F \oplus (\Delta^p_{M2})^F \oplus
(\Delta^p_{M3})^F,$$ where, for $i=1,2,3$, $(\Delta^p_{Mi})^F$ is
the Friedrichs extension of the restriction of $\Delta^p_M$ to
$C^{\infty}_c(\Lambda^p((c,+\infty)\times N))\cap {\mathcal
L}_i((c,+\infty)\times N).$\par Since the orthogonal sum is
finite, for $1\leq p \leq n-1$, $$ \sigma_{\rm
ess}((\Delta^p_M)^F)= \bigcup_{i=1}^3 \sigma_{\rm
ess}((\Delta^p_{Mi})^F),$$ $$
\sigma_p((\Delta^p_M)^F)=\bigcup_{i=1}^3 \sigma_p
((\Delta^p_{Mi})^F).$$ For $p=0$ (resp. $p=n$), any $\omega \in
L^2_p((c,+\infty)\times N))$ can be written as
$\omega=\omega_{1\delta}$ (resp. $\omega= \omega_{2d} \wedge dt$),
where $\omega_{1\delta}$ (resp. $\omega_{2d}$) is a coclosed
(resp. closed) $0$-form (resp. $(n-1)$-form) parametrized by $t$
on $N$. Hence $L^2_0((c,+\infty)\times N)= {\mathcal
L}_1((c,+\infty)\times N)$ (resp. $L^2_n((c,+\infty)\times N)=
{\mathcal L}_2((c,+\infty)\times N)$), and $(\Delta^p_M)^F=
(\Delta^p_{M1})^F$ (resp. $(\Delta^p_M)^F=(\Delta^p_{M2})^F$).
\par Hence, for any $p\in [0,n]$, in order to determine the
spectrum of $(\Delta^p_M)^F$ it suffices to study the spectral
properties of $(\Delta^p_{Mi})^F$, $i=1,2,3$.
\par To this purpose, let us introduce further decompositions.
First of all, we decompose $\omega_{1\delta}$ according to an
orthonormal basis $\left\{\tau_{1k}\right\}_{k \in \mathbb N}$ of
coclosed $p$-eigenforms of $\Delta^p_N$; this yields
\begin{equation}\label{1delta} \omega_{1\delta}=\oplus_k
h_k(t) \tau_{1k},\end{equation} where $h_k(t)\tau_{1k}\in
L^2_p((c,+\infty)\times N)$ for every $k \in \mathbb N$. The
sum (\ref{1delta}) is orthogonal in $L^2_p((c,+\infty)\times N)$, thanks
to
(\ref{exp}). We will call $p$-form of type I any $p$-form
$\omega\in L^2_p((c,+\infty)\times N)$ such that $$\omega= h(t)
\tau_1,$$ where $\tau_1$ is a coclosed normalized $p$-eigenform,
corresponding to some eigenvalue $\lambda$ of $\Delta^p_N$. For
every $k \in \mathbb N$, let us denote by $\lambda_k^p$ the
eigenvalue associated to $\tau_{1k}$. Since for every $k\in
\mathbb N$
\begin{multline}\label{h1}\Delta^p_{M1}(h(t)\tau_{1k})=
\frac{\lambda^p_k}{g(t)} h(t)\tau_{1k} \\-
f(t)^{-\frac{1}{2}}g(t)^{\frac{-n+1+2p}{2}}\frac{\partial}{\partial
t}\left(f(t)^{-\frac{1}{2}}g(t)^{\frac{n-1-2p}{2}}\frac{\partial
h}{\partial t}\right) \tau_{1k}, \end{multline} $\Delta^p_{M1}$ is
invariant under decomposition (\ref{1delta}), and, since, if
$\omega=h(t)\tau_{1k}$, $$\|\omega\|^2_{L^2_p((c,+\infty)\times
N)}=
\int_c^{\infty}g(s)^{\frac{n-2p-1}{2}}f(s)^{\frac{1}{2}}h(s)^2\,ds,$$
$(\Delta^p_{M1})^F$ is unitarily equivalent to the direct sum over $k
\in \mathbb N$ of the Friedrichs extensions
$(\Delta_{1\lambda^p_k})^F$ of the operators
$$\Delta_{1\lambda^p_k}: C^{\infty}_c (c,+\infty) \longrightarrow
L^2((c,+\infty), g^{\frac{n-2p-1}{2}}f^{\frac{1}{2}})$$
\begin{multline} \Delta_{1\lambda^p_k}h=\left\{
\frac{\lambda^p_k}{g(t)} h(t)  -
f(t)^{-\frac{1}{2}}g(t)^{\frac{-n+1+2p}{2}}\frac{\partial}{\partial
t}\left(f(t)^{-\frac{1}{2}}g(t)^{\frac{n-1-2p}{2}}\right)
\right\}.\end{multline} If we introduce the tranformation
\begin{equation}\label{trasf1}w(t)=
h(t)f(t)^{\frac{1}{4}}g(t)^{\frac{n-2p-1}{4}},
\end{equation}
a direct (but lengthy) computation shows that $(\Delta^p_{M1})^F$ is
unitarily equivalent to the direct sum, over $k\in \mathbb N$, of
the Friedrichs extensions $(D_{1\lambda^p_k})^F$ of the operators
$$D_{1\lambda^p_k}: C^{\infty}_c(c,+\infty)\longrightarrow
L^2(c,+\infty)$$ given by
\begin{multline}\label{w1} D_{1\lambda^p_k}w = -
\frac{\partial}{\partial t}\left( \frac{1}{f} \frac{\partial
w}{\partial t} \right)+ \left\{ -\frac{7}{16} \frac{1}{f^3}
\left(\frac{\partial f}{\partial t}\right)^2+ \frac{1}{4}
\frac{1}{f^2} \frac{\partial^2 f}{\partial t^2}+\right.\\ -
\frac{1}{2} \frac{1}{f^2}\frac{\partial f}{\partial
t}\frac{(n-1-2p)}{4} \frac{1}{g}\frac{\partial g}{\partial t} +
\frac{1}{f} \frac{(n-2p-1)}{4}
\frac{(n-2p-5)}{4}\frac{1}{g^2}\left(\frac{\partial g}{\partial
t}\right)^2+\\ \left. +
\frac{1}{f}\frac{(n-2p-1)}{4}\frac{1}{g}\frac{\partial^2 g
}{\partial t^2} + \frac{\lambda^p_k}{g}\right\} w.
\end{multline}
Analogously, we decompose $\omega_{2d}$ according to an
orthonormal basis of closed $(p-1)$-eigenforms $\left\{\tau_{2k}
\right\}_{k \in \mathbb N}$ of $\Delta^{p-1}_N$:
\begin{equation}\label{2d} \omega_{2d}\wedge dt=\oplus_k
h_k (t) \tau_{2k}\wedge dt.\end{equation} We will call $p$-form of
type II any $p$-form $\omega\in L^2_p((c,+\infty)\times N)$ such
that $$\omega= h(t) \tau_2 \wedge dt,$$ where
$\tau_2$ is a
coclosed normalized $(p-1)$-eigenform, corresponding to some
eigenvalue $\lambda$ of $\Delta^{p-1}_N$. For every $k \in \mathbb N$
$$\Delta^p_{M2}(h(t)\tau_{2k}\wedge dt )= (\Delta_{2\lambda^{p-1}_k}
h )\tau_{2k}\wedge dt,$$ where
\begin{multline}\label{h2}
\Delta_{2\lambda^{p-1}_k} h = \frac{\lambda^{p-1}_k}{g(t)} h(t)-
\\ \frac{\partial}{\partial t}\left\{
f(t)^{-\frac{1}{2}}g(t)^{\frac{-n-1+2p}{2}}
\frac{\partial}{\partial t}\left( f(t)^{-\frac{1}{2}}
g(t)^{\frac{n+1-2p}{2}}h(t) \right) \right\}. \end{multline} Here,
again, for every $k \in \mathbb N$ we denote by $\lambda_k^{p-1}$
the eigenvalue of $\Delta^{p-1}_N$ corresponding to the eigenform
$\tau_{2k}$. Since if $\omega= h(t) \tau_{2k}\wedge dt$
$$\|\omega\|^2_{L^2_p((c,+\infty)\times N)}=
\int_c^{\infty}g(s)^{\frac{n-2p+1}{2}}f(s)^{-\frac{1}{2}}h(s)^2\,ds,$$
introducing the transformation
\begin{equation}\label{trasf2} w(t)= h(t)
f(t)^{-\frac{1}{4}}g(t)^{\frac{n+1-2p}{4}},
\end{equation}
we find that $(\Delta^p_{M2})^F$ is unitarily equivalent to the direct
sum, over $k \in \mathbb N$, of the Friedrichs extensions $(D_{2
\lambda^{p-1}_k})^F$ of the operators $$D_{2 \lambda^{p-1}_k}:
C^{\infty}_c(c,+\infty) \longrightarrow L^2(c,+\infty) $$
\begin{multline}\label{w2} D_{2\lambda^{p-1}_k}w = -
\frac{\partial}{\partial t}\left( \frac{1}{f} \frac{\partial
w}{\partial t} \right)+  \left\{ -\frac{7}{16} \frac{1}{f^3}
\left(\frac{\partial f}{\partial t}\right)^2+
\frac{1}{4} \frac{1}{f^2} \frac{\partial^2 f}{\partial t^2}+\right.\\ -
\frac{1}{2} \frac{1}{f^2}\frac{\partial f}{\partial
t}\frac{(n-1+2p)}{4} \frac{1}{g}\frac{\partial g}{\partial t} +
\frac{1}{f} \frac{(n-2p+1)}{4}
\frac{(n-2p+5)}{4}\frac{1}{g^2}\left(\frac{\partial g}{\partial
t}\right)^2+ \\ \left.
+\frac{1}{f}\frac{(-n+2p-1)}{4}\frac{1}{g}\frac{\partial^2
g}{\partial t^2} + \frac{\lambda^{p-1}_k}{g}\right\} w.
\end{multline}
Decompose now $\omega_{2\delta}$ with respect to an
orthonormal basis of coclosed $(p-1)$-eigenforms $\left\{
\tau_{3k}\right\}_{k \in \mathbb N}$ of $\Delta^{p-1}_N$; for every
$k \in \mathbb N$, denote by $\lambda_k^{p-1}$ the eigenvalue
corresponding to the eigenform $\tau_{3k}$. \par Then
$\left\{\frac{1}{\sqrt{\lambda^{p-1}_k}}d_N \tau_{3k} \right\}_{k
\in \mathbb N}$ is an orthonormal basis of closed eigenforms
 on $N$, and we get the following decomposition for
$\omega_{1d} \oplus \omega_{2\delta}\wedge dt$:
\begin{multline}\label{3}\omega_{1d}\oplus \omega_{2\delta}\wedge dt =
\oplus_k (\frac{1}{\sqrt{\lambda^{p-1}_k}} h_{1k} d^{p-1}_N \tau_{3k}
\oplus (-1)^p h_{2k}\tau_{3k}\wedge dt).
\end{multline}
We will call $p$-form of type III any $p$-form $\omega \in
L^2_p((c,+\infty)\times N)$ such that $$\omega=
\frac{1}{\sqrt{\lambda}}h_1(t)d^{p-1}_N \tau_3 \oplus_M (-1)^p h_2(t)
\tau_3 \wedge dt ,$$ where $\tau_3$ is a normalized coclosed
$(p-1)$-eigenform of $\Delta^{p-1}_N$, corresponding to the eigenvalue
$\lambda$. A direct computation shows that
\begin{multline}\label{h1h2}
\Delta^p_{M3}\left(\frac{1}{\sqrt{\lambda}} h_1(t)d^{p-1}_N \tau_3
\oplus_M (-1)^p h_2(t) \tau_3 \wedge dt\right)=\\ =\left(\Delta_{1\lambda}
h_1 + \frac{1}{f(t)}\frac{1}{g(t)} \frac{\partial g}{\partial
t}\sqrt{\lambda} h_2 \right) \left( \frac{1}{\sqrt{\lambda}} d^{p-1}_N
\tau_3 \right)\\ \oplus \left(\Delta_{2\lambda} h_2 + \frac{1}{g^2(t)}
\frac{\partial g}{\partial t} \sqrt{\lambda} h_1 \right)
\left((-1)^p \tau_3 \wedge dt \right);
\end{multline}
moreover, if $\omega=\frac{1}{\sqrt{\lambda}} h_1(t)d^{p-1}_N \tau_3
\oplus_M (-1)^p h_2(t) \tau_3 \wedge dt$, then
$$\|\omega\|^2_{L^2_p((c,+\infty)\times N)}= \int_c^{+\infty}
g(s)^{\frac{n-2p-1}{2}}f(s)^{\frac{1}{2}} h_1(s)^2 \,ds $$ $$+
\int_c^{+\infty} g(s)^{\frac{n+1-2p}{2}}f(s)^{-\frac{1}{2}}
h_2(s)^2\,ds. $$ Hence, introducing the transformation
\begin{equation}\label{trasf12}\begin{array}{ll}
w_1(t)=&g^{\frac{n-2p-1}{4}}(t)f^{\frac{1}{4}}(t)h_1(t) \\
w_2(t)=& g^{\frac{n-2p+1}{4}}(t)f^{-\frac{1}{4}}(t) h_2(t),
\end{array}
\end{equation}
we find that $(\Delta^p_{M3})^F$ is unitarily equivalent to the direct
sum, over $k \in \mathbb N$, of the Friedrichs extensions
$(D_{3\lambda^{p-1}_k})^F$ of the operators
$$D_{3\lambda^{p-1}_k}:C^{\infty}_c(c,+\infty) \oplus
C^{\infty}_c(c,+\infty) \longrightarrow L^2(c,+\infty)\oplus
L^2(c,+\infty) $$
\begin{multline}\label{w1w2}D_{3\lambda^{p-1}_k}(w_1\oplus w_2)=
\left(D_{1\lambda^{p-1}_k}w_1 +
g(t)^{-\frac{3}{2}}f(t)^{-\frac{1}{2}}\frac{\partial g}{\partial
t} \sqrt{\lambda^{p-1}_k}w_2 \right) \oplus \\ \oplus
\left(D_{2\lambda^{p-1}_k} w_2 +
g(t)^{-\frac{3}{2}}f(t)^{-\frac{1}{2}}\frac{\partial g}{\partial
t} \sqrt{\lambda^{p-1}_k}w_1 \right).
\end{multline}
Hence the study of the spectrum of $\Delta^p_M$ can be reduced to
the analysis of the spectral properties of the selfadjoint
operators $(D_{1\lambda_k^p})^F$, $(D_{2\lambda_k^{p-1}})^F$ and
$(D_{3\lambda_k^{p-1}})^F$. Since the Hodge $*$ operator maps
$p$-forms of type I isometrically onto $(n-p)$-forms of type II,
it suffices to consider the cases $i=1$ and $i=3$.\par
We have that, for $i=1,3$
$$\sigma_{\rm
ess}((\Delta^p_{Mi})^F)\supset \bigcup_k \sigma_{\rm
ess}((D_{i\lambda_k})^F).$$
However, since the direct sums in (\ref{1delta}) and (\ref{3}) have an
infinite number of summands, we cannot conclude that $$\sigma_{\rm
ess}((\Delta^p_{Mi})^F)= \bigcup_k \sigma_{\rm
ess}((D_{i\lambda_k})^F)$$
without additional information on the isolated eigenvalues of
$(D_{i\lambda_k})^F$.

\section{The case $a=-1$}
For $a=-1$, the operators $D_{1\lambda^p_k}$,
$D_{2\lambda^{p-1}_k}$, $D_{3\lambda^{p-1}_k}$ are given by $$
D_{1\lambda^p_k}w= -\frac{\partial^2 w}{\partial t^2}+\left[
\left(\frac{n-2p-1}{2} \right)^2b^2 +\lambda^p_k
e^{2bt}\right]w,$$ $$D_{2\lambda^{p-1}_k}w= - \frac{\partial^2
w}{\partial t^2} + \left[\left(\frac{n-2p+1}{2}\right)^2 b^2
+\lambda^{p-1}_ke^{2bt}\right]w,$$
$$D_{3\lambda^{p-1}_k}(w_1\oplus w_2)=$$ $$=  -\frac{\partial^2
w_1}{\partial t^2}+\left[ \left(\frac{n-2p-1}{2} \right)^2b^2
+\lambda^{p-1}_k e^{2bt}\right]w_1 -2b e^{bt}
\sqrt{\lambda^{p-1}_k} w_2 \oplus$$ $$ \oplus  -\frac{\partial^2
w_2}{\partial t^2}+\left[ \left(\frac{n-2p+1}{2} \right)^2b^2
+\lambda^{p-1}_k e^{2bt}\right]w_2 -2b
e^{bt}\sqrt{\lambda^{p-1}_k}w_1.$$ The behaviour of the potential
part depends strongly on the sign of $b$, hence we will
investigate separately the case $b<0$, $b=0$, $b>0$.

\subsection{The case $b<0$}
\par \bigskip \noindent
For $a=-1$, $b<0$, the situation is similar to the asymptotically
hyperbolic case treated in \cite{Antoci}. The volume of $M$ is
infinite, and the sectional curvatures
$K_M(\frac{\partial}{\partial t},\frac{\partial}{\partial
\theta_i})$, $K_M(\frac{\partial}{\partial
\theta_i},\frac{\partial}{\partial \theta_j})$ tend to $-b^2$ as
$t\rightarrow +\infty$. \par First of all, we will compute the
essential spectrum of $(\Delta^p_{M1})^{F}$. To this purpose, let
us observe that for every $k \in \mathbb N$, ${\mathcal
D}((D_{1\lambda^p_k})^F)\subset W^{1,2}_0(c,+\infty)$. Indeed, if
$\left\{w_n\right\}$ is a Cauchy sequence in
$C^{\infty}_c(c,+\infty)$ with respect to the graph norm, that is,
if for every $\epsilon>0$ there exists $\bar n$ such that for
every $n,m>\bar n$
$$\|w_n-w_m\|_{L^2(c,+\infty)}+\|D_{1\lambda^p_k}(w_n-w_m)
\|_{L^2(c,+\infty)}<\epsilon ,$$ then for every $n,m>\bar n$
$$\int_c^{+\infty}\left(1+\left(\frac{n-2p-1}{2}\right)^2b^2+\lambda^p_k
e^{2bs}\right)(w_n-w_m)^2\, ds+$$ $$+\int_c^{+\infty}
\left(\frac{\partial}{\partial s}(w_n-w_m)\right)^2 \,ds
<\epsilon. $$ Hence, $\left\{ w_n\right\}$ is a Cauchy sequence in
$W^{1,2}(c,+\infty)$. As a consequence, if $w\in {\mathcal
D}((D_{1\lambda^p_k})^F)$, $w \in W^{1,2}_0(c,+\infty)$. \par We
are now in position to prove our first result.
\begin{lem}\label{-1<01l} Let $a=-1$, $b<0$. For $0\leq p\leq n-1$, for every
$k\in \mathbb N$, $$\sigma_{\rm ess}((D_{1\lambda^p_k})^F)=
\left[ \left( \frac{n-2p-1}{2}\right)^2 b^2, +\infty \right).$$
\end{lem}
\begin{proof}
Let us consider the Friedrichs extension $(D_{10})^F$ of the
operator  $$ D_{10}:C^{\infty}_c (c,+\infty)\longrightarrow
L^2(c,+\infty)$$ $$ D_{10}w= -\frac{\partial^2w}{\partial
t^2}+\left(\frac{n-2p-1}{2}\right)^2 b^2 w.$$ It is not difficult
to check that, for every $k \in \mathbb N$,
$(D_{1\lambda^p_k})^F-(D_{10})^F$ is a relatively compact
perturbation of $(D_{10})^F$. \par Indeed, $${\mathcal
D}((D_{10})^F)\subseteq {\mathcal
D}((D_{1\lambda^p_k})^F-(D_{10})^F)={\mathcal D}((D_{10})^F),$$
since $X_{D_{1\lambda^p_k}}= X_{D_{10}}$ and ${\mathcal
D}(D_{1\lambda^p_k}^*)= {\mathcal D}(D_{10}^*)$. \par Moreover, if
$\left\{w_n\right\}\subset {\mathcal D}((D_{01})^F)$ satisfies the
condition $$
\|w_n\|_{L^2(c,+\infty)}+\|(D_{01})^Fw_n\|_{L^2(c,+\infty)}\leq
C,$$ $\left\{w_n\right\}$ is bounded in $W^{1,2}(c,+\infty)$;
hence it is bounded in $L^{\infty}(c,+\infty)$ and in $W^{1,2}(K)$
for every compact subset $K\subset (c,+\infty)$. Thus, since for
$b<0$ $e^{2bt}\in L^2(c,+\infty)\cap L^{\infty}(c,+\infty)$, for
every $\tilde c >c$ and for every $n,m\in \mathbb N$ $$
\int_c^{+\infty}(\lambda^p_k)^2 e^{4bs} (w_n-w_m)^2\,ds \leq$$ $$
\leq C_1 \int_c^{\tilde c} (w_n-w_m)^2 \,ds +C_2 \int_{\tilde
c}^{+\infty}(\lambda^p_k)^2 e^{4bs}\,ds,$$ with $C_1,C_2$ not
depending on $\tilde c$. Let us consider a sequence $\left\{c_h
\right\}\subset (c,+\infty)$ such that $c_h \rightarrow +\infty$
as $h\rightarrow +\infty$ and for every $h \in \mathbb N$
$$\int_{c_h}^{+\infty}(\lambda^p_k)^2 e^{4bs}\,ds <\frac{1}{h} .$$
For $h=1$, thanks to the Rellich-Kondrachov theorem, there exists
a subsequence $\left\{w_{n(1)}\right\}\subseteq \left\{w_n
\right\}$ which converges in $L^2(c,\tilde c_1)$. Hence, for every
$\eta>0$ there exists $\bar n(1)$ such that for every $n,m>\bar
n(1)$ $$\int_c^{+\infty}(\lambda^p_k)^2
e^{4bs}(w_{n(1)}-w_{m(1)})^2\,ds
< \eta +1 .$$ Analogously, for $h=2$ there exists a subsequence
$\left\{w_{n(2)}\right\}\subseteq \left\{ w_{n(1)} \right\}$ such
that for every $\eta >0$ there exists $\bar n(2)$ for which
$$n,m>\bar n(2) \Longrightarrow \int_c^{+\infty}(\lambda^p_k)^2
e^{4bs}(w_{n(2)}-w_{m(2)})^2 \,ds
< \eta +\frac{1}{2}.$$ Iterating this argument, for every $h\in \mathbb N$
we find
a subsequence \\ $\left\{w_{n(h)}\right\} \subseteq \left\{
w_{n(h-1)} \right\}$ such that for every $\eta >0$ there exists
$\bar n(h)$ for which $$n,m>\bar n(h) \Longrightarrow
\int_c^{\infty}(\lambda^p_k)^2 e^{4bs}(w_{n(h)}-w_{m(h)})^2\,ds<
\eta +\frac{1}{h}.$$ Through a Cantor diagonal process, we can
then find a subsequence \\ $\left\{w_{n_l}\right\}\subset
\left\{w_n\right\}$,
 such that for every $\eta>0$ and every $h\in
\mathbb N$,
if $l,m>\bar l= \bar n(h)$ $$\|((D_{1\lambda^p_k})^F-(D_{10})^F)
(w_{n_l}-w_{n_m)})\|^2_{L^2(c,+\infty)}=$$
$$=\int_c^{+\infty}(\lambda^p_k)^2 e^{4bs}(w_{n_l}-w_{n_m})\,ds <
\eta +\frac{1}{h}. $$ Hence, $(D_{1\lambda^p_k})^F- (D_{10})^F$ is a
relatively compact perturbation of $(D_{10})^F$. \par Since
$\sigma_{\rm ess}(D_{10}^F)= [(\frac{n-2p-1}{2})^2b^2,+\infty)$,
the conclusion follows.
\end{proof}
As a consequence, $\left[ \left( \frac{n-2p-1}{2}\right)^2 b^2,
+\infty \right) \subseteq \sigma_{\rm ess}((\Delta^p_{M1})^F)$. On the
other hand,
\begin{lem}
Let $a=-1$, $b<0$. Then, if $\mu < \left( \frac{n-2p-1}{2}
\right)^2 b^2$, $\mu \notin \sigma_{\rm ess}((\Delta^p_{M1})^F)$.
\end{lem}
\begin{proof} If $\mu \in \sigma_{\rm ess }((\Delta^p_{M1})^F)$, there
exists a Weyl sequence $\left\{\omega_k\right\}\subset {\mathcal
D}((\Delta^p_{M1})^F)$ for $\mu$. $\left\{\omega_k\right\}$ has no
convergent subsequence and we have
$$ \langle \omega_k,\omega_k \rangle_{L^2_p((c,+\infty)\times N)} \leq
C,$$
$$ \lim_{k\rightarrow +\infty}((\Delta^p_{M1})^F\omega_k -\mu
\omega_k)=0.$$
Moreover, we can
suppose that
$$\omega_k= h_{\lambda^p_k}\tau_{\lambda^p_k} ,$$
where, for every $k\in \mathbb N$, $\tau_{\lambda^p_k}$
 is a normalized coclosed
$p$-eigenform of $\Delta^p_N$ corresponding to the eigenvalue
$\lambda^p_k$
and $\lambda^p_k \rightarrow +\infty$ as $k\rightarrow +\infty$.\par
Hence, there exists a bounded sequence $\left\{w_k\right\}$ in
$L^2(c,+\infty)$
such that,
for every $k \in \mathbb N$, $w_k \in {\mathcal D}((D_{1\lambda^p_k})^F)$,
and such that
$$\lim_{k\rightarrow +\infty}\|(D_{1\lambda^p_k})^F w_k
-\mu w_k\|_{L^2(c,+\infty)}=0, $$
from which we cannot extract any convergent subsequence.\par
Then
$$\langle (D_{1\lambda^p_k})^F w_k- \mu
w_k,w_k\rangle_{L^2(c,+\infty)}\longrightarrow 0 $$
as $k\rightarrow +\infty$ and since for every $k\in \mathbb N$
$${\mathcal D}((D_{1\lambda^p_k})^F)\subset W^{1,2}_0(c,+\infty),$$ then
$$ \int_c^{+\infty} \left(\frac{\partial w_k}{\partial s} \right)^2\,ds
+\int_c^{+\infty}\lambda^p_k
e^{2bs} w_k^2 \,ds +$$ $$+ \int_c^{+\infty}\left( \left(\frac{n-2p-1}{2}
\right)^2 b^2 -\mu \right) w_k^2 \, ds \longrightarrow 0 $$
as $k \rightarrow +\infty$. Thus
$$\int_c^{+\infty}w_k^2\,ds \longrightarrow 0 $$
as $k\rightarrow +\infty$, and we get a contradiction. Hence, if $\mu <
\left(\frac{n-2p-1}{2}\right)^2b^2$, $\mu$ cannot belong to the essential
spectrum of
$(\Delta^p_{M1})^F$.
\end{proof}
As a consequence, we have the following
\begin{prop}\label{-1<0M1} Let $a=-1$, $b<0$. For $0\leq p\leq n-1$,$$\sigma_{\rm
ess}((\Delta^p_{M1})^F) = \left[\left(\frac{n-2p-1}{2}\right)^2
b^2,+\infty \right).$$ \end{prop} By duality,
\begin{prop}\label{-1<0M2} Let $a=-1$, $b<0$. For $1\leq p\leq n$, $$\sigma_{\rm
ess}((\Delta^p_{M2})^F)= \left[
\left(\frac{n-2p+1}{2}\right)^2b^2,+\infty\right).$$ \end{prop}We
still have to investigate the spectrum of $(\Delta^p_{M3})^F$. First
of all, we compute the essential spectrum of
$(D_{3\lambda_k^{p-1}})^F$ for every $k \in \mathbb N$.
\begin{lem} Let $a=-1$, $b<0$. For $1\leq p \leq n-1$ and for every $k \in \mathbb N$,
\begin{multline} \sigma_{\rm ess}((D_{3\lambda^{p-1}_k})^F)=\\ =\left[\min
\left\{\left(\frac{n-2p-1}{2}\right)^2 b^2
,\left(\frac{n-2p+1}{2}\right)^2b^2\right\},+\infty \right).\end{multline}
\end{lem}
\begin{proof}
Let us consider the Friedrichs extension $(D_{30})^F$ of the
operator $$ D_{30}:C^{\infty}_c(c,+\infty) \oplus C^{\infty}_c
(c,+\infty) \longrightarrow L^2(c,+\infty)\oplus L^2(c,+\infty)$$
$$ D_{30}(w_1\oplus w_2)= -\frac{\partial^2 w_1}{\partial
t^2}+\left( \frac{n-2p-1}{2}\right)^2b^2 w_1 \oplus$$ $$ \oplus -
\frac{\partial^2 w_2}{\partial t^2} +
\left(\frac{n-2p+1}{2}\right)^2b^2 w_2.$$ Since the essential
spectrum of $(D_{30})^F$ is equal to $$ \left[ \min \left\{
\left(\frac{n-2p-1}{2}\right)^2 b^2,
\left(\frac{n-2p+1}{2}\right)^2b^2 \right\}, +\infty \right),$$ it
suffices to show that $(D_{3\lambda^{p-1}_k})^F-(D_{30})^F$ is a
relatively compact perturbation of $(D_{30})^F$.\par Now,
$${\mathcal D}((D_{30})^F) \subseteq {\mathcal
D}((D_{3\lambda^{p-1}_k})^F-(D_{30})^F)={\mathcal D}((D_{30})^F),
$$ since $X_{D_{30}}=X_{D_{3\lambda^{p-1}_k}}$ and ${\mathcal
D}(D^*_{30})={\mathcal D}(D^*_{3\lambda^{p-1}_k})$. \par Moreover,
if the sequence $\left\{ w_{1n}\oplus w_{2n}\right\}$ is such
that, for every $n\in \mathbb N$, $$ \|(w_{1n}\oplus
w_{2n})\|_{L^2(c,+\infty)\oplus
L^2(c,+\infty)}+\|D_{30}^F(w_{1n}\oplus
w_{2n})\|_{L^2(c,+\infty)\oplus L^2(c,+\infty)}\leq C,$$ then
$\left\{w_{1n}\right\}$ and $\left\{w_{2n}\right\}$ are bounded in
$W^{1,2}(c,+\infty)$; hence, they are bounded also in
$L^{\infty}(c,+\infty)$ and in $W^{1,2}(K)$ for any compact subset
$K\subset (c,+\infty)$.\par Since $b<0$, $e^{2bt}\in
L^2(c,+\infty) \cap L^{\infty}(c,+\infty)$ and $e^{bt}\in
L^2(c,+\infty)\cap L^{\infty}(c,+\infty)$; thus, arguing as in the
proof of Lemma \ref{-1<01l}, we can construct a subsequence
$\left\{w_{1n_h}\oplus w_{2n_h}\right\}$ such that for every
$\delta >0$ there exists $\bar h$ such that whenever $h,l>\bar h$
$$\|((D_{3\lambda_k^{p-1}})^F-(D_{30})^F)((w_{1n_h}-w_{1n_l})
\oplus(w_{2n_h}-w_{2n_l})) \|_{L^2(c,+\infty)\oplus
L^2(c,+\infty)}=$$ $$=\|\lambda_k^{p-1} e^{2bt}
(w_{1n_h}-w_{1n_l})-2b
\sqrt{\lambda_k^{p-1}}e^{bt}(w_{2n_h}-w_{2n_l})\|_{L^2(c,+\infty)}+$$
$$ +\|\lambda_k^{p-1} e^{2bt}(w_{2n_h}-w_{2n_l}) -2b
\sqrt{\lambda_k^{p-1}}e^{bt}(w_{1n_h}-w_{1n_l})\|_{L^2(c,+\infty)}<\delta.$$
Hence, $(D_{3\lambda^{p-1}_k})^F-(D_{30})^F$ is a relatively
compact perturbation of $(D_{30})^F$, and $\sigma_{\rm
ess}((D_{3\lambda^{p-1}_k})^F)=\sigma_{\rm ess}((D_{30})^F)$.
 \end{proof}
On the other hand,
\begin{lem}\label{difficile} Let $a=-1$, $b<0$. If $$0 <\mu <
\min\left\{\left(\frac{n-2p-1}{2}\right)^2 b^2 ,
\left(\frac{n-2p+1}{2}\right)^2 b^2 \right\},$$ then $\mu \notin
\sigma_{\rm ess}((\Delta^p_{M3})^F)$.
\end{lem}
\begin{proof}
Let
$$0<\mu<\min\left\{\left(\frac{n-2p-1}{2}\right)^2b^2,
\left(\frac{n-2p+1}{2}\right)^2b^2
\right\};$$ then $\mu \in \sigma_{\rm ess}((\Delta_{M3}^p)^F)$ if and only
if
there exist a sequence $\left\{\mu_k\right\}$ of eigenvalues of
$(\Delta^p_{M3})^F$
and a corresponding sequence of normalized, mutually orthogonal,
eigenforms
$\left\{\Phi_k\right\}$ of $(\Delta_{M3}^p)^F$ such that for
every
$k\in \mathbb N$
$$ (\Delta^p_{M3})^F \Phi_k - \mu_k \Phi_k =0$$ and
$$ \mu_k \longrightarrow \mu \quad \quad \hbox{as $k\rightarrow
+\infty$}.$$
In view of the weak Kodaira decomposition, replacing
$\left\{\Phi_k\right\}$ by a subsequence (again denoted by the same
symbol for shortness), we can suppose that either
$\delta^p_M \Phi_k =0$ for every $k\in \mathbb N$, or
$d^p_M \Phi_k=0$ for every $k\in \mathbb N$. \par \noindent
a) In the first case, since $\mu \not= 0$, for every $k\in \mathbb N$
 $d^p_M \Phi_k \not=0$; we have $\|d^p_M
\Phi_k\|_{L^2_{p+1}((c,+\infty)\times N)}<C$, and $$
\Delta^{p+1}_M d^p_M \Phi_k - \mu \,d^p_M \Phi_k \longrightarrow
0$$ as $k \rightarrow +\infty$ because $$ \|\Delta^{p+1}_M
d^p_M\Phi_k - \mu \,d^p_M\Phi_k \|_{L^2_{p+1}((c,+\infty)\times
N)}\leq$$ $$\|\Delta^{p+1}_M d^p_M \Phi_k - \mu_k \,d^p_M
\Phi_k\|_{L^2_{p+1}((c,+\infty)\times
N)}+|\mu_k-\mu|\|d^p_M\Phi_k\|_{L^2_{p+1}((c,+\infty)\times N)},$$
where $$ \Delta^{p+1}_M d^p_M \Phi_k - \mu_k \,d^p_M \Phi_k=0$$
for every $k\in \mathbb N$, and $\mu_k \rightarrow \mu$ as
$k\rightarrow +\infty$. Moreover, $\left\{d^p_M \Phi_k\right\}$
has no convergent subsequences, since for $i\not= j$ $\langle
d^p_M \Phi_i,d^p_M \Phi_j\rangle_{L^2_{p+1}((c,+\infty)\times N)}=
0$ and $\|d^p_M\Phi_k\|^2=\mu_k > \mu -\epsilon$ for $k$ big
enough.\par Now, let $\psi \in C^{\infty}_c (M)$ be such that
$$0\leq \psi(x)\leq 1 \quad \hbox{for every $x\in M$},$$
$$\psi(x)=1 \quad \hbox{for $x \in M\setminus ((c,+\infty)\times
N)$}$$ and $$\psi(x)=0 \quad \hbox{for $x \in (2c,+\infty)\times
N$}.$$ If we set, for every $k\in \mathbb N$, $$
\tilde{\omega}_k:= (1-\psi)( d^p_M \Phi_{2k+1}- d^p_M
\Phi_{2k}),$$ arguing as in the proof of Satz 3.1 in
\cite{Eichhorn} we find that $\left\{\tilde{\omega}_k\right\}$ is
a Weyl sequence for $\mu$ for the operator $(\Delta^{p+1}_M)^F$.
Since an explicit computation shows that, for every $k\in \mathbb
N$, $d^p_M \Phi_k$ is a $(p+1)$-form of type II and
$(1-\psi)d^p_M\Phi_k \in {\mathcal D}((\Delta^{p+1}_{M2})^F)$, we
can argue that $\mu$ lies in the essential spectrum  of
$(\Delta^{p+1}_{M2})^F$. But this is not possible because by
assumption $$\mu < \left(\frac{n-2p-1}{2}\right)^2b^2=
\left(\frac{n-2(p+1)+1}{2}\right)^2b^2.$$ \noindent b) Consider
now the second case. Since $\mu>0$, for every $k\in \mathbb N$
$\delta^p_M \Phi_k \not=0$; following the same argument as in part
a), if we set for every $k\in \mathbb N$ $$\tilde{\omega}_k:=
(1-\psi)(\delta^p_M \Phi_{2k+1}- \delta^p_M \Phi_{2k}), $$ where
$\psi(x)\in C^{\infty}_c(M)$ is chosen as in part a), we find that
$\left\{\tilde{\omega}_k\right\}$ is a Weyl sequence for $\mu$ for
$(\Delta^{p-1}_M)^F$. Since an explicit computation shows that,
for every $k\in \mathbb N$, $\delta^p_M \Phi_k$ is a $(p-1)$-form
of type I, we conclude that $\mu$ lies in the essential spectrum
of $(\Delta^{p-1}_{M1})^F$.
 But this is not possible since by assumption
$$\mu <\left(\frac{n-2p+1}{2}\right)^2b^2=
\left(\frac{n-2(p-1)-1}{2}\right)^2b^2
.$$
 \end{proof}
In the proof of Lemma \ref{difficile}, it is essential that
$\left\{\mu_k\right\}$ is a sequence of strictly positive real
numbers. Hence, by this technique we cannot determine whether $0$
is an isolated eigenvalue of infinite multiplicity of
$(\Delta_{M3}^p)^F$ or not.
\par As a consequence, all we can state about the essential
spectrum of $(\Delta_{M3}^p)^F$ is:
\begin{prop}\label{-1<0M3} Let $a=-1$, $b<0$. Then, for $1\leq p \leq
n-1$,
\begin{multline}\sigma_{\rm ess}((\Delta^p_{M3})^F)\setminus
\left\{0\right\}=\\= \left[ \min
\left\{\left(\frac{n-2p-1}{2}\right)^2b^2,\left(\frac{n-2p+1}{2}\right)^2
b^2 \right\},+\infty\right).\end{multline} \end{prop} Combining
the results of Proposition \ref{-1<0M1}, Proposition \ref{-1<0M2}
and Proposition \ref{-1<0M3}, finally we can state the following
\begin{thm}\label{a=-1b<0} Let $a=-1$, $b<0$, $0\leq p\leq n$. Then, if $p
\not= \frac{n\pm 1}{2}$ \begin{multline}\sigma_{\rm
ess}(\Delta^p_M)\setminus \left\{0\right\} =\\ =\left[
\min\left\{\left(\frac{n-2p-1}{2}\right)^2
b^2,\left(\frac{n-2p+1}{2}\right)^2 b^2 \right\},+\infty
\right),\end{multline} whilst if $p=\frac{n+1}{2}$ or
$p=\frac{n-1}{2}$ $$\sigma_{\rm ess}(\Delta_M)=[0,+\infty).$$
\end{thm}
\bigskip \bigskip
\subsection{The case $b=0$}\par \noindent
For $a=-1$, $b=0$, the operators $D_{1\lambda^p_k}$,
$D_{2\lambda^{p-1}_k}$, $D_{3\lambda^{p-1}_k}$ are simply $$
D_{1\lambda^p_k}w = - \frac{\partial^2 w}{\partial t^2}+
\lambda^p_k w,$$ $$ D_{2\lambda^{p-1}_k}w= - \frac{\partial^2
w}{\partial t^2} +\lambda^{p-1}_k w,$$ $$ D_{3\lambda^{p-1}_k}(w_1
\oplus w_2)= - \frac{\partial^2 w_1}{\partial t^2} +
\lambda^{p-1}_k w_1 \oplus - \frac{\partial^2 w_2}{\partial
t^2}+\lambda^{p-1}_k w_2;$$ hence the essential spectrum of
$D_{1\lambda^p_k}^F$ (resp. $D_{3\lambda^{p-1}_k}^F$) is equal to
$[\lambda^p_k,+\infty)$ (resp. $[\lambda^{p-1}_k,+\infty)$).
Moreover, an explicit computation shows that, for every $k\in
\mathbb N$, $(D_{1\lambda_k^p})^F$ (resp.
$(D_{3\lambda_k^{p-1}})^F$) has no eigenvalues. As a consequence,
we recover the well-known result for cylindrical ends (see e.g.
\cite{Melrose}):
\begin{thm}\label{a=-1b=0} Let $a=-1$, $b=0$. Then, for $0\leq p \leq n$, the
essential spectrum of $\Delta^p_M$ is given by $$ \sigma_{\rm
ess}(\Delta^p_M)=\bigcup_k ([\lambda^p_k,+\infty)\cup
[\lambda^{p-1}_k,+\infty))= [\overline\lambda, +\infty),$$ where
$\overline \lambda= \min_k \left\{
\lambda^p_k,\lambda^{p-1}_k\right\}$.
\end{thm}
In particular, if the $p$-th or the $(p-1)$-th Betti number of $N$
does not vanish, $\overline \lambda=0$, otherwise $\overline \lambda
>0$. \bigskip \bigskip

\subsection{The case $b>0$}
As in the previous cases, we begin with the spectral analysis of
$D_{1\lambda^p_k}^F$ for every $k\in \mathbb N$:
\begin{lem}\label{-1>01k} Let $a=-1$, $b>0$. Then, for every $k\in \mathbb N$,
if $\lambda^p_k=0$, $\sigma_{\rm ess}((D_{1\lambda^p_k})^F)=
\left[ \left( \frac{n-2p-1}{2} \right)^2 b^2, +\infty \right)$;
if, on the contrary, $\lambda^p_k >0$, $\sigma_{\rm
ess}((D_{1\lambda^p_k})^F)= \emptyset$.
\end{lem}
\begin{proof}
The first assertion is obvious; as for the other, it is well-known (see
e.g.
\cite{Muller-Pfeiffer}, Thm. 3.13) that the spectrum of any selfadjoint
extension of
an operator of type $-\frac{\partial^2}{\partial t^2} + V(t)$ acting on
$C^{\infty}_c(c,+\infty)$ is purely discrete if, and only if, for every
$h \in (0,1)$
$$\lim_{t\rightarrow +\infty}\int_t^{t+h} V(s)\,ds = +\infty  .$$
Since for $\lambda^p_k >0$
$$ \int_t^{t+h}\left[\left(\frac{n-2p-1}{2}\right)^2 b^2 +\lambda^p_k
e^{2bs} \right]\,ds =$$
$$= \left( \frac{n-2p-1}{2} \right)^2b^2 h +\frac{\lambda^p_k}{2b}
e^{2bt}(e^{2bh}-1)\longrightarrow +\infty$$
as $t\rightarrow +\infty$, the conclusion follows.
\end{proof}
Now, an easy computation shows that for every $k \in \mathbb N$
\begin{equation}\label{stima-1>0} \langle (D_{1\lambda^p_k})^F
w,w\rangle_{L^2(c,+\infty)} \geq
\left(\left(\frac{n-2p-1}{2}\right)^2b^2
+\lambda^p_k \right) \langle w,w
\rangle_{L^2(c,+\infty)}.\end{equation}
As a consequence, $\sigma_p((\Delta_{M1}^p)^F)$ have no cluster points and
every
eigenvalue of $(\Delta^p_{M1})^F$ has finite multiplicity. Hence:
\begin{prop}\label{betti1} Let $a=-1$, $b>0$. For $0\leq p \leq n-1$, if the
$p$-th Betti number of $N$ vanishes, $\sigma_{\rm
ess}((\Delta^p_{M1})^F)=
\emptyset$;
if, on the contrary, the $p$-th Betti number of $N$ is
different from zero, $\sigma_{\rm ess}((\Delta^p_{M1})^F)= \left[
\left(\frac{n-2p-1}{2}\right)^2 b^2, +\infty \right).$ \end{prop}
By duality,
\begin{prop}\label{betti2} Let $a=-1$, $b>0$. For $1\leq p \leq n$, if the
$(p-1)$-th Betti number of $N$ vanishes, $\sigma_{\rm
ess}((\Delta^p_{M2})^F)= \emptyset$; if, on the contrary, the
$(p-1)$-th Betti number of $N$ is different from zero,
$\sigma_{\rm ess}((\Delta^p_{M2})^F)= \left[
\left(\frac{n-2p+1}{2}\right)^2b^2, +\infty \right)$. \end{prop}
We shall now investigate the spectrum of $(\Delta^p_{M3})^F$. As a
first step, we will compute the essential spectrum of
$(D_{3\lambda_k^{p-1}})^F$ for every $k\in \mathbb N$. To this
purpose, we need a preliminary Lemma:
\begin{lem}\label{fittizio} For every $K \in \Real$ and for every
$\lambda \geq 0$, the
essential spectrum of the Friedrichs extension $D^F$ of the
operator
$$ D: C^{\infty}_c(c,+\infty)\oplus C^{\infty}_c(c,+\infty)\longrightarrow
L^2(c,+\infty)\oplus L^2(c,+\infty)$$
$$ {\mathcal D}(w_1\oplus w_2)= - \frac{\partial^2w_1}{\partial t^2}+
\left[
K+\lambda e^{2bt}\right] w_1 -2b \sqrt{\lambda}e^{bt}w_2 \oplus$$
$$ \oplus - \frac{\partial^2 w_2}{\partial t^2}+\left[K+ \lambda e^{2bt}
\right]w_2 -2b \sqrt{\lambda} e^{bt}w_1$$
is empty.
\end{lem}
\begin{proof} First of all, through an argument similar to that of Satz
3.1 in
\cite{Eichhorn}, it can be shown that the essential
spectrum of $D^F$ does not depend on the choice of the first endpoint $c$
of $(c,+\infty)$. Hence, given $K$ and $\mu\geq 0$ we can
suppose that for every $t>c$
\begin{equation}\label{posy}K- \mu + \lambda e^{2bt} \pm 2b e^{bt}
\sqrt{\lambda}> C>0.\end{equation}
Consider now the closed subspaces ${\mathcal V}_1$, ${\mathcal V}_2$ of
$L^2(c,+\infty)\oplus
L^2(c,+\infty)$ defined as
$${\mathcal V}_1:=\left\{w_1\oplus w_2\,|\,w_1=w_2 \right\} ,$$
$${\mathcal V}_2:= \left\{w_1\oplus w_2 \,|\, w_1=-w_2 \right\}.$$
They are orthogonal in
$L^2(c,+\infty)\oplus L^2(c,+\infty)$ and any $w_1 \oplus w_2$ can be
written as
$$w_1 \oplus w_2 = \left(\frac{w_1 + w_2}{2} \oplus \frac{w_1 +w_2}{2}
\right) + \left(\frac{w_1 -w_2}{2} \oplus
\frac{w_2 -w_1}{2} \right). $$
Hence,
\begin{equation}\label{decc}L^2(c,+\infty) \oplus L^2(c,+\infty)=
{\mathcal V}_1\oplus{\mathcal V}_2. \end{equation} Moreover, an
explicit computation shows that the operator $D$ is invariant
under decomposition (\ref{decc}). As a consequence, $D^F$ splits
as $$D^F=(D_{|{\mathcal V}_1})^F \oplus (D_{|{\mathcal V}_2})^F,
$$ where, for $i=1,2$, $D_{|{\mathcal V}_i}$ is the restriction of
$D$ to $$\left(C^{\infty}_c(c,+\infty) \oplus
C^{\infty}_c(c,+\infty)\right)\cap {\mathcal V}_i;$$ thus,
$$\sigma_{\rm ess}(D^F)=\sigma_{\rm ess}(D_{|{\mathcal
V}_1}^F)\cup \sigma_{\rm ess}(D_{|{\mathcal V}_2}^F) .$$ Let us
begin with $D_{|{\mathcal V}_1}^F$. If $\mu \in \sigma_{\rm
ess}(D_{|{\mathcal V}_1}^F)$, there exists a Weyl sequence for
$\mu$, that is a sequence $\left\{ w_n \right\} \subset {\mathcal
D}(D_{|{\mathcal V}_1}^F)$ such that $\|w_n\|_{L^2(c,+\infty)}\leq
C$ and $(D - \mu)(w_n \oplus w_n) \longrightarrow 0 $ as $n
\rightarrow +\infty$, with no convergent subsequences. Then
$$\langle (D-\mu)(w_n \oplus w_n),(w_n \oplus
w_n)\rangle_{L^2(c,+\infty)\oplus L^2(c,+\infty)}\longrightarrow 0
$$ as $n\rightarrow +\infty$. Thus $$\int_c^{+\infty}\left(
\frac{\partial w_n}{\partial s}\right)^2\,ds +\int_c^{+\infty}
\left[K-\mu +\lambda
e^{2bs}-2b\sqrt{\lambda}e^{bs}\right]w_n^2\,ds \longrightarrow 0$$
as $n\rightarrow +\infty$. Since (\ref{posy}) implies that
$$\int_c^{+\infty}w_n^2\,ds < C^{-1}\int_c^{+\infty}\left[K-\mu
+\lambda e^{2bs}- 2b \sqrt{\lambda}e^{bs}\right]w_n^2\,ds
 ,$$ we find that
$\|w_n \|_{L^2(c,+\infty)}\rightarrow 0$ as $n\rightarrow +\infty$. Hence,
$\sigma_{\rm ess}(D_{|{\mathcal V}_1}^F)= \emptyset$.\par
Suppose now that $\mu \in \sigma_{\rm ess}(D_{|{\mathcal V
}_2}^F)$.
If $\left\{ w_n\right\}$ is a Weyl sequence for $\mu$ for $D_{|{\mathcal
V}_2}^F$, we have $$\langle (D - \mu)(w_n \oplus -w_n),(w_n \oplus - w_n)
\rangle_{L^2(c,+\infty)\oplus L^2(c,+\infty)}\longrightarrow 0 $$
as $n\rightarrow +\infty$. Hence
$$\int_c^{+\infty}\left( \frac{\partial w_n}{\partial s}\right)^2\,ds +
\int_c^{+\infty}\left[ K- \mu +\lambda e^{2bs}+ 2b \sqrt{\lambda}e^{bs}
\right]w_n^2 \,ds \longrightarrow 0 $$
as $n \rightarrow +\infty$. In view of (\ref{posy}),
$\|w_n\|_{L^2(c,+\infty)} \rightarrow 0$ as $n\rightarrow +\infty$. But a
Weyl sequence cannot converge. Hence, $\sigma_{\rm
ess}(D_{|{\mathcal V}_2}^F)=\emptyset$.
\end{proof}
We can now compute the essential spectrum of $(D_{3\lambda_k^{p-1}})^F$:
\begin{lem}\label{-1>03k} Let $a=-1$, $b>0$. Then, for $1\leq p\leq n-1$ and for any $k \in \mathbb
N$,
$$\sigma_{\rm
ess}((D_{3\lambda^{p-1}_k})^F)=
\emptyset.$$
\end{lem}
\begin{proof} Let us consider the Friedrichs extension $(D_{30})^F$ of the
operator
$$ D_{30}:C^{\infty}_c(c,+\infty) \oplus C^{\infty}_c(c,+\infty)
\longrightarrow L^2(c,+\infty)\oplus L^2(c,+\infty)$$
\begin{multline} D_{30}(w_1 \oplus w_2)= - \frac{\partial^2 w_1}{\partial
t^2}
+ \tilde
K w_1 +\lambda^{p-1}_k e^{2bt} w_1   -2b \sqrt{\lambda^{p-1}_k}
e^{bt}w_2
\oplus \\
 \oplus - \frac{\partial^2 w_2}{\partial t^2}+ \tilde K w_2 +
\lambda^{p-1}_k e^{2bt} w_2 -2b \sqrt{\lambda^{p-1}_k} e^{bt}
w_1,\end{multline} where $$\tilde K= \max\left\{
\left(\frac{n-2p-1}{2}\right)^2b^2,
\left(\frac{n-2p+1}{2}\right)^2b^2 \right\}.$$ From the previous
Lemma, we know that $\sigma_{\rm ess}((D_{30})^F)= \emptyset$. We
will show that $((D_{3\lambda^{p-1}_k})^F- (D_{30})^F)$ is a
relatively compact perturbation of $(D_{30})^F$.\par First of all,
since a straightforward computation shows that
$X_{D_{30}}\subseteq X_{D_{3\lambda^{p-1}_k}}$ and ${\mathcal
D}(D_{30}^*)\subseteq {\mathcal D}(D_{3\lambda^{p-1}_k}^*)$, then
$${\mathcal D}((D_{30})^F)\subseteq {\mathcal
D}((D_{3\lambda^{p-1}_k})^F-(D_{30})^F) ={\mathcal
D}((D_{3\lambda^{p-1}_k})^F).$$ \par \noindent We still have to
show that, given a sequence $\left\{w_{1n}\oplus w_{2n}\right\}$
in ${\mathcal D}((D_{30})^F)$ such that $$ \|w_{1n}\oplus
w_{2n}\|_{L^2(c,+\infty)}^2 + \|((D_{30})^F)(w_{1n}\oplus
w_{2n})\|_{L^2(c,+\infty)}^2\leq C$$ there exists a subsequence
$\left\{w_{1n_k}\oplus w_{2n_k}\right\}$ such that
$$\left\{((D_{3\lambda^{p-1}_k})^F- (D_{30})^F)(w_{1n_k}\oplus
w_{2n_k})\right\}$$ converges.\par Now, the fact that
$\|(D_{30})^F(w_{1n}\oplus w_{2n})\|^2_{L^2(c,+\infty)}\leq C$ is
equivalent to the inequalities
\begin{equation}\label{uno} \|- \frac{\partial^2 w_{1n}}{\partial t^2}
+\left[
\tilde K + \lambda_k^{p-1} e^{2bt}\right]w_{1n}-2b
e^{bt}\sqrt{\lambda_k^{p-1}}w_{2n}\|_{L^2(c,+\infty)}\leq C,\end{equation}
\begin{equation}\label{due} \|- \frac{\partial^2 w_{2n}}{\partial t^2}
+\left[
\tilde K + \lambda_k^{p-1} e^{2bt}\right]w_{2n}-2b
e^{bt}\sqrt{\lambda_k^{p-1}}w_{1n}\|_{L^2(c,+\infty)}\leq C,\end{equation}
which in turn imply
\begin{multline}\label{piu}\|-\frac{\partial^2 (w_{1n}+w_{2n})}{\partial
t^2}+\left[\tilde K + \lambda_k^{p-1}
e^{2bt}\right](w_{1n}+w_{2n})+ \\  -2b
e^{bt}\sqrt{\lambda_k^{p-1}}(w_{2n}+w_{1n})\|_{L^2(c,+\infty)}\leq C,
\end{multline}
\begin{multline}\label{meno}\|-\frac{\partial^2 (w_{1n}-w_{2n})}{\partial
t^2}+\left[\tilde K +  \lambda_k^{p-1}
e^{2bt}\right](w_{1n}-w_{2n})+ \\ -2b
e^{bt}\sqrt{\lambda_k^{p-1}}(w_{2n}-w_{1n})\|_{L^2(c,+\infty)}\leq C.
\end{multline}
By taking the inner product with $(w_{1n}+w_{2n})$ and
$(w_{1n}-w_{2n})$ respectively, (\ref{piu}) and (\ref{meno})
yield:
\begin{multline}\label{piu1}\int_c^{+\infty}\left(\frac{\partial
(w_{1n}+w_{2n})}{\partial s}\right)^2\,ds +\\
+\int_c^{+\infty}\left(\tilde K + \lambda_k^{p-1}
e^{2bs}- 2b \sqrt{\lambda_k^{p-1}}e^{bs}\right)(w_{1n}+w_{2n})^2\,ds \leq
C,
\end{multline}
\begin{multline}\label{meno1}
\int_c^{+\infty}\left(\frac{\partial
(w_{1n}-w_{2n})}{\partial s}\right)^2\,ds +\\+
\int_c^{+\infty}\left(\tilde K + \lambda_k^{p-1}
e^{2bs}+ 2b \sqrt{\lambda_k^{p-1}}e^{bs}\right)(w_{1n}-w_{2n})^2\,ds \leq
C.
\end{multline}
Now, there exists $\bar c$ such that for every
strictly positive $\lambda^{p-1}_k\in
\sigma_p(\Delta_N)$ and for every $t>\bar c$
$$\tilde K + \lambda^{p-1}_ke^{2bt}\pm 2b \sqrt{\lambda^{p-1}_k}e^{bt}>
\tilde
C e^{\frac{3}{2}bt} ,$$
with $\tilde C>0$. Since in view of the results of \cite{Eichhorn} we can
assume $c>\bar c$, we obtain from (\ref{piu1}) and (\ref{meno1}) that
$$\int_c^{+\infty}e^{\frac{3}{2}bs}(w_{1n}\pm w_{2n})^2\,ds \leq C .$$
Hence, for $i=1,2$,
$$\int_c^{+\infty}e^{\frac{3}{2}bs}(w_{in})^2\,ds \leq C .$$
Moreover, $\left\{ w_{1n}\right\}$ and $\left\{w_{2n}\right\}$ are bounded
in $W^{1,2}(c,+\infty)$; thus, they are bounded also in
$L^{\infty}(c,+\infty)$ and in $W^{1,2}(K)$ for every compact subset
$K\subset (c,+\infty)$. As a consequence, if we set
$$K_1(n,p)=\left(\frac{n-2p-1}{2}\right)^2 b^2,$$
$$K_2(n,p)= \left( \frac{n-2p+1}{2}\right)^2 b^2,$$
for every $n,m \in \mathbb N$
$$\|((D_{3\lambda^{p-1}_k})^F-(D_{30})^F)((w_{1n}-w_{1m})
\oplus(w_{2n}-w_{2m}))\|^2_{L^2(c,+\infty)}
=$$
$$=\sum_{i=1,2}\int_c^{+\infty}(K_i(p,n)-\tilde K)^2(w_{in}-w_{im})^2\,ds
\leq $$
$$ \leq \sum_{i=1,2}(K_i(p,n)-\tilde
K)^2\|e^{-\frac{3}{4}bt}(w_{in}-w_{im})\|_{L^2}\|e^{\frac{3}{4}bt}
(w_{in}-w_{im})\|_{L^2}\leq$$
$$\leq C
\sum_{i=1,2}\|e^{-\frac{3}{4}bt}(w_{in}-w_{im})\|_{L^2(c,+\infty)} .$$
Since $$e^{-\frac{3}{4}bt}\in L^2(c,+\infty)\cap L^{\infty}(c,+\infty),$$
following the argument of  Lemma \ref{-1<01l}, we can extract from
$\left\{
w_{1n}\oplus w_{2n} \right\}$  a
subsequence $\left\{w_{1n_k}\oplus w_{2n_k}\right\}$ such
that $$\left\{((D_{3\lambda_k^{p_1}})^F-(D_{30})^F)(w_{1n_k}\oplus
w_{2n_k})\right\}$$ converges. \par Hence, $\sigma_{\rm
ess}(D_{3\lambda_k^{p-1}})^F)= \sigma_{\rm
ess}(D_{30}^F)= \emptyset$.
\end{proof}
We still have to investigate the cluster points of
$\sigma_p((\Delta^p_{M3})^F)$, which could be additional points in
the essential spectrum of $\Delta_M^p$. The following Lemma holds:
\begin{lem}\label{ficicchia} Let $a=-1$, $b>0$. Let $0<\mu \in
\overline{\sigma_p((\Delta^p_{M3}))^F}$.
\begin{enumerate}
\item If the $p$-th and the $(p-1)$-th Betti numbers of $N$ both vanish,
$\mu$
is an isolated eigenvalue of finite multiplicity;  \\
\item if the $p$-th Betti number of $N$ vanishes whilst the $(p-1)$-th
Betti number of $N$ is different from zero, and if $\mu$ is not an
isolated eigenvalue of finite multiplicity, then $\mu \geq
\left(\frac{n-2p+1}{2}\right)^2b^2$;\\
\item if the $p$-th Betti number of $N$ is different from zero whilst
the $(p-1)$-th Betti number of $N$ vanishes, and if $\mu$ is not an
isolated eigenvalue of finite multiplicity, then $\mu \geq
\left(\frac{n-2p-1}{2}\right)^2b^2$; \\
\item if the $p$-th and the $(p-1)$-th Betti numbers of $N$ are both
different from zero, and if $\mu$ is not an isolated eigenvalue of finite
multiplicity, then $\mu \geq \min
\left\{\left(\frac{n-2p-1}{2}\right)^2b^2,\left(\frac{n-2p+1}{2}\right)^2
b^2 \right\}$. \end{enumerate} \end{lem}
\begin{proof}
If $\mu$ is an eigenvalue of infinite multiplicity of
$(\Delta^p_{M3})^F$ or is a cluster point of
$\sigma_p((\Delta^p_{M3})^F)$, there exist a sequence
$\left\{\mu_k\right\}$ of eigenvalues of $(\Delta^p_{M3})^F$ and a
corresponding sequence of normalized, mutually orthogonal
eigenforms $\left\{\Phi_k\right\}$ such that  for every $k\in
\mathbb N$ $$ \Delta^p_M \Phi_k - \mu_k \Phi_k =0$$ and $$ \mu_k
\longrightarrow \mu \quad \quad \hbox{as $ k \rightarrow
+\infty$.}$$ In view of the weak Kodaira decomposition, replacing
$\left\{\Phi_k\right\}$ by a subsequence (again denoted by the
same symbol for shortness) we can suppose that either $\delta^p_M
\Phi_k=0$ for every $k \in \mathbb N$, or $d^p_M \Phi_k =0$ for
every $k\in \mathbb N$. Following the argument of Lemma
\ref{difficile}, we see that, in the first case, $\mu$ lies in the
essential spectrum of the operator $(\Delta^{p+1}_{M2})^F$,
whilst, in the second case, $\mu$ lies in the essential spectrum
of the operator $(\Delta^{p-1}_{M1})^F$.\par \noindent a) Consider
the first case. If the $p$-th Betti number of $N$ vanishes, we get
a contradiction with Proposition \ref{betti2}; if, on the
contrary, the $p$-th Betti number of $N$ is different from zero,
Proposition \ref{betti2} implies that $$\mu \geq
\left(\frac{n-2(p+1)+1}{2}\right)^2b^2=\left(
\frac{n-2p-1}{2}\right)^2b^2 .$$ \noindent b) Consider the second
case. If the $(p-1)$-th Betti number of $N$ vanishes, we get a
contradiction with Proposition \ref{betti1}; if, on the contrary,
the $(p-1)$-th Betti number of $N$ is different from zero, by
Proposition \ref{betti2} we have that $$\mu \geq
\left(\frac{n-2(p-1)-1}{2}\right)^2b^2
=\left(\frac{n-2p+1}{2}\right)^2b^2.$$ Combining these facts we
complete the proof.\end{proof} As in the proof of Lemma
\ref{difficile}, it is essential that $\mu
>0$.  As a consequence, all we can say about the essential spectrum of
$(\Delta^p_{M3})^F$ is:
\begin{prop}\label{-1>03M} Let $a=-1$, $b>0$, $1\leq p \leq n-1$.
\begin{enumerate}
\item
If both the $p$-th Betti number and the $(p-1)$-th Betti number of $N$
vanish, $\sigma_{\rm ess}((\Delta^p_{M3})^F)\setminus \left\{0\right\} =
\emptyset.$\\
\item If the $p$-th Betti number of $N$ is different from zero  and the
$(p-1)$-th Betti number of $N$
vanishes, $$ \sigma_{\rm ess}((\Delta^p_{M3})^F)\setminus
\left\{0\right\}\subseteq \left[
\left(\frac{n-2p-1}{2}\right)^2b^2,+\infty\right).$$\\
\item If the $p$-th Betti number of $N$ vanishes and the $(p-1)$-th Betti
number of $N$ is different from zero, $$ \sigma_{\rm
ess}((\Delta^p_{M3})^F)\setminus \left\{ 0\right\} \subseteq \left[\left(
\frac{n-2p+1}{2}\right)^2 b^2, +\infty \right).$$\\
\item If both the $p$-th and the $(p-1)$-th Betti numbers of $N$ are
different from zero,
$$\sigma_{\rm ess}((\Delta^p_{M3})^F)\setminus \left\{ 0 \right\}
\subseteq$$
$$\subseteq
\left[\min\left\{\left(\frac{n-2p-1}{2}\right)^2b^2,\left(\frac{n-2p+1}
{2}\right)^2b^2\right\},+\infty\right).$$
\end{enumerate}
\end{prop}
Combining Proposition \ref{betti1}, Proposition \ref{betti2} and
Proposition \ref{-1>03M}, we finally get the following Theorem,
which shows an interesting link between the cohomology of the
boundary $N$ and the essential spectrum of $\Delta_M^p$ (with $0$
excluded):
\begin{thm} Let $a=-1$, $b>0$, $0\leq p \leq n$.
\begin{enumerate}
\item
If both the $p$-th and the $(p-1)$-th Betti numbers of $N$ vanish,
$$ \sigma_{\rm ess}(\Delta^p_M)\setminus \left\{0\right\} = \emptyset.$$\\
\item
If the $p$-th Betti number of $N$ is different from zero  and the
$(p-1)$-th Betti number of $N$
vanishes, $$ \sigma_{\rm ess}(\Delta^p_M)\setminus \left\{0\right\}=
\left[ \left(\frac{n-2p-1}{2}\right)^2b^2,+\infty\right).$$\\
\item If the $p$-th Betti number of $N$ vanishes and the $(p-1)$-th Betti
number of $N$ is different from zero, $$ \sigma_{\rm
ess}(\Delta^p_M)\setminus \left\{ 0\right\} = \left[\left(
\frac{n-2p+1}{2}\right)^2 b^2, +\infty \right).$$\\
\item If both the $p$-th and the $(p-1)$-th Betti numbers of $N$ are
different from zero,
$$\sigma_{\rm ess}(\Delta^p_M)\setminus \left\{ 0 \right\} =$$$$=
\left[\min\left\{\left(\frac{n-2p-1}{2}\right)^2b^2,\left(\frac{n-2p+1}
{2}\right)^2b^2\right\},
+\infty\right). $$ \end{enumerate}
\end{thm}

\section{The case $a<-1$}
Let us perform the change of variables $$ (0,+\infty)
\longrightarrow \left(\frac{1}{|a+1|},+\infty\right) $$ $$r(t):=
\frac{e^{-(a+1)t}}{|a+1|}. $$ The Riemannian metric in the new
coordinate system $(r,\theta)$ on \\ $(\tilde c, +\infty)\times
N$, where $\tilde c = r(c)$, takes the form $${ds'}^2=
dr^2+|a+1|^{\frac{2b}{a+1}}\,r^{\frac{2b}{a+1}}\,d\theta_N^2.$$ If
we apply the decomposition techniques in the new coordinate
system, we find the operators $D_{1\lambda^p_k}$,
$D_{2\lambda^{p-1}_k}$, $D_{3\lambda^{p-1}_k}$, defined on the
smooth functions with compact support in $(\tilde c,+\infty)$:
 $$D_{1\lambda^p_k}w= - \frac{\partial^2
w}{\partial r^2}+ K_1(n,p) r^{-2}w+
\lambda^p_k|a+1|^{\frac{-2b}{a+1}}\,r^{\frac{-2b}{a+1}} w,$$
$$D_{2\lambda^{p-1}_k}w= -\frac{\partial^2 w}{\partial r^2}+
K_2(n,p) r^{-2}w +
\lambda^{p-1}_k|a+1|^{\frac{-2b}{a+1}}\,r^{\frac{-2b}{a+1}}w ,$$
$$D_{3\lambda^{p-1}_k}(w_1\oplus w_2)= \left( - \frac{\partial^2
w_1}{\partial r^2}+ K_1(n,p) r^{-2}w_1+
\lambda^{p-1}_k|a+1|^{\frac{-2b}{a+1}}r^{\frac{-2b}{a+1}} w_1
+\right.$$ $$\left.
+\sqrt{\lambda^{p-1}_k}|a+1|^{-\frac{b}{a+1}}r^{-\frac{b}{a+1}-1}w_2
\right) \oplus $$ $$\oplus \left( - \frac{\partial^2 w_2}{\partial
r^2}+ K_2(n,p) r^{-2}w_2+
\lambda^p_k|a+1|^{\frac{-2b}{a+1}}r^{\frac{-2b}{a+1}} w_2
+\right.$$ $$\left.
+\sqrt{\lambda^{p-1}_k}|a+1|^{-\frac{b}{a+1}}r^{-\frac{b}{a+1}-1}w_1\right)
,$$ where $$K_1(n,p)= \left(\frac{n-2p-1}{2}\right)^2
\frac{b^2}{(a+1)^2} +\frac{n-2p-1}{2}\frac{b}{|a+1|} ,$$
$$K_2(n,p)= \left(\frac{n-2p+1}{2}\right)^2 \frac{b^2}{(a+1)^2}
+\frac{n-2p+1}{2}\frac{b}{|a+1|} .$$ Since the potential terms
containing $p$ and $n$ tend to zero as $r\rightarrow +\infty$, we
can presume that for $a=-1$ the bottom of the essential spectrum
of $\Delta^p_M$ will not depend on the relationships between the
dimension $n$ of $M$ and the degree $p$. The asymptotic behaviour
of the potential is again strongly determined by the sign of $b$.
Hence also for $a<-1$ we will consider separately the cases $b<0$,
$b=0$, $b>0$.
\subsection{The case $b<0$} \par \noindent
For $a<-1$, $b<0$, all the potential terms in the operators
$D_{1\lambda_k^p}$, $D_{2\lambda_k^{p-1}}$, $D_{3\lambda_k^{p-1}}$
tend to zero as $r\rightarrow +\infty$; hence the following result
is not surprising:
\begin{thm}\label{a<-1b<0} Let $a<-1$, $b<0$. Then, for every $p\in [0,n]$, $$
\sigma_{\rm
ess}(\Delta^p_M)=[0,+\infty).$$
\end{thm}
\begin{proof} Since $\sigma(\Delta_M^p)\subseteq [0,+\infty]$, it suffices
to show that, for every $k\in \mathbb N$,
$\sigma_{\rm
ess}(D_{1\lambda^p_k})= [0,+\infty)$.\par
To this purpose, let us consider the Friedrichs extension $(D_{10})^F$ of
the Laplacian $-\frac{\partial^2}{\partial r^2}$ on $C^{\infty}_c(\tilde
 c, +\infty)$. We will show that $(D_{1\lambda^p_k})^F- (D_{10})^F$ is a
relatively compact perturbation of $(D_{10})^F$. \par First of
all, since $X_{D_{1\lambda^p_k}}= X_{D_{10}}$ and ${\mathcal
D}(D_{1\lambda^p_k}^*)= {\mathcal D}(D_{10}^*)$, we have that
${\mathcal D}((D_{10})^F)\subseteq {\mathcal
D}((D_{1\lambda^p_k})^F - (D_{10})^F)= {\mathcal
D}((D_{10})^F)$.\par Now, $(D_{1\lambda^p_k})^F-(D_{10})^F$ is
given by $$(D_{1\lambda^p_k}^F-D_{10}^F)w=( K_1(n,p) r^{-2}+
\lambda^p_k |a+1|^{-\frac{2b}{a+1}}\,r^{-\frac{2|b|}{|a+1|}})\,w,
$$ where $$K_1(n,p)r^{-2}\in L^2(\tilde c,+\infty)\cap
L^{\infty}(\tilde c,+\infty)$$ and $$\lambda^p_k \,
|a+1|^{-\frac{2|b|}{|a+1|}} r^{-\frac{2|b|}{|a+1|}} \in
L^{\infty}_\eps (\tilde c,+\infty)$$ for $0< \eps \leq
\frac{|b|}{|a+1|}$ since $$(1+r^2)^{\eps} r^{-\frac{2|b|}{|a+1|}}
\leq 2^{\eps}r^{2(\eps-\frac{|b|}{|a+1|})}.$$ If
$\left\{w_n\right\}\subset {\mathcal D}((D_{01})^F)$ is such that
$$\|w_n\|_{L^2(\tilde c,+\infty)}+ \|(D_{01})^Fw_n\|_{L^2(\tilde
c,+\infty)}\leq C,$$ then $$\|w_n\|_{L^2(\tilde
c,+\infty)}^2+\|\frac{\partial w_n}{\partial r}\|_{L^2(\tilde
c,+\infty)}^2\leq C,$$ hence $\left\{w_n\right\}$ is a bounded
sequence in $W^{1,2}(\tilde c,+\infty)$, in $L^{\infty}(\tilde
c,+\infty)$ and in $W^{1,2}(K)$ for every compact subset $K
\subset (\tilde c,+\infty)$.\par We have
$$\|((D_{1\lambda^p_k})^F-(D_{10})^F)(w_n-w_m)\|_{L^2(\tilde
c,+\infty)} \leq $$ $$\leq C_1 \|r^{-2}(w_n-w_m)\|_{L^2(\tilde
c,+\infty)}+\lambda^p_k
C_2\|r^{-\frac{2|b|}{|a+1|}}(w_n-w_m)\|_{L^2(\tilde c,+\infty)}
.$$ As for the first term, since $r^{-2}\in L^2(\tilde c,+\infty)
\cap L^{\infty}(\tilde c,+\infty)$, following the same argument as
in Lemma \ref{-1<01l} we find a subsequence, again denoted by
$\left\{w_n\right\}$, such that for every $\eta >0$ there exists
$\bar n$ for which $$ n,m>\bar n \quad \Longrightarrow C_1
\|r^{-2}(w_n-w_m) \|_{L^2(\tilde c,+\infty)}<\frac{\eta}{2}.$$ As
for the second term, we cannot apply the same argument because
$r^{-\frac{2|b|}{|a+1|}}$ might not belong to $L^2(\tilde
c,+\infty)$. Nevertheless, since $r^{-\frac{2|b|}{|a+1|}}\leq
C(1+r^2)^{-\eps}$ for $0<\eps
<\frac{|b|}{|a+1|}$, we have that, for every $d >\tilde c$, $$
\int_{\tilde c}^{+\infty}x^{-\frac{4|b|}{|a+1|}}(w_n-w_m)^2\,dx
\leq C \int_{\tilde c}^{+\infty}(1+x^2)^{-2\eps}(w_n-w_m^2)^2\,dx
=$$ $$ = C\int_{\tilde c}^{d}(1+x^2)^{-2\eps}(w_n-w_m)^2\,dx +
C\int_d^{+\infty}(1+x^2)^{-2\eps}(w_n-w_m)^2\,dx\leq$$ $$ \leq
C\|w_n-w_m\|_{L^2(\tilde c,d)}^2 + C\|w_n-w_m\|_{L^{\infty}(\tilde
c,+\infty)} \frac{1}{(1+d^2)^{2\eps}} \leq $$ $$\leq
C\|w_n-w_m\|_{L^2(\tilde c,d)}^2 + \tilde C
\frac{1}{(1+d^2)^{2\eps}}. $$ for some positive constant $\tilde
C$. Let $\left\{\tilde c_h\right\} \subset (\tilde c,+\infty)$ be
a sequence such that $\tilde c_h \rightarrow +\infty$ as
$h\rightarrow +\infty$ and, for every $h\in \mathbb N$,
$$\frac{\tilde C}{(1+\tilde c_h^2)^{2\eps}}< \frac{1}{h}.$$ Then,
again through an argument similar to that of Lemma \ref{-1<01l},
we can extract a subsequence, again denoted by
$\left\{w_n\right\}$ for shortness, such that, for every $\eta>0$
and for every $h\in \mathbb N$, there exists $\bar n(h)$ for which
$$n,m> \bar n(h) \quad \Longrightarrow \quad \int_{\tilde
c}^{+\infty}x^{-\frac{4|b|}{|a+1|}}(w_n-w_m)^2 \,dx <
\frac{\eta}{2} +\frac{1}{h}.$$ Hence, we have found a subsequence
$\left\{w_{n}\right\}$ such that for every $\eta >0$ and for every
$h>0$, there exists $\tilde n = \max\left\{ \bar n, \bar n(h)
\right\}$ for which $$n,m >\tilde n \quad \Longrightarrow \quad
\|((D_{1\lambda^p_k})^F-(D_{10})^F) (w_n-w_m)\|_{L^2(\tilde
c,+\infty)} \leq \eta + \frac{1}{h} .$$ As a consequence,
$(D_{1\lambda^p_k})^F-(D_{10})^F$ is a relatively compact
perturbation of $(D_{10})^F$.  \end{proof}

\subsection{The case $b=0$} For $a<-1$, $b=0$, the operators
$D_{1\lambda^p_k}$, $D_{2\lambda^{p-1}_k}$, $D_{3\lambda^{p-1}_k}$
are simply $$ D_{1\lambda^p_k}w = - \frac{\partial^2 w}{\partial
r^2}+ \lambda^p_k w,$$ $$ D_{2\lambda^{p-1}_k}w= -
\frac{\partial^2 w}{\partial r^2} +\lambda^{p-1}_kw,$$ $$
D_{3\lambda^{p-1}_k}(w_1 \oplus w_2)= \left( - \frac{\partial^2
w_1}{\partial r^2} + \lambda^{p-1}_k w_1 \right) \oplus \left( -
\frac{\partial^2 w_2}{\partial r^2}+\lambda^{p-1}_k w_2\right).$$
Hence, arguing as in the case $a=-1$, $b=0$ we get the following
result:
\begin{thm}\label{a<-1b=0} Let $a<-1$, $b=0$. Then, for $0\leq p \leq n$,
the essential spectrum of $\Delta^p_M$ is given by $$ \sigma_{\rm
ess}(\Delta^p_M)=\bigcup_k
([\lambda^p_k,+\infty)\cup[\lambda^{p-1}_k,+\infty))=
[\overline \lambda, +\infty),$$ where $\overline \lambda=
\min_k\left\{\lambda^p_k,\lambda^{p-1}_k\right\}$.
\end{thm} In
particular, if the $p$-th or the $(p-1)$-th Betti number of $N$ does
not vanish, $\overline \lambda=0$. Otherwise, $\overline \lambda >0$.

\subsection{ The case $b>0$}
Once more, we begin with the spectral analysis of
$(D_{1\lambda^p_k})^F$ for every $k\in \mathbb N$:
\begin{lem} Let $a<-1$, $b>0$. For any $k\in \mathbb N$, if $\lambda^p_k=0$, $\sigma_{\rm
ess}((D_{1\lambda^p_k})^F)=
\left[0, +\infty \right)$; if, on the contrary, $\lambda^p_k >0$,
 $\sigma_{\rm ess}((D_{1\lambda^p_k})^F)= \emptyset$.
\end{lem}
\begin{proof}
If $\lambda^p_k=0$, $$D_{1\lambda^p_k}= -\frac{\partial^2}{\partial r^2}+
K_1(n,p)r^{-2},$$ and since $$r^{-2}\in L^2(\tilde c,+\infty) \cap
L^{\infty}(\tilde c,+\infty),$$ the essential spectrum of
$(D_{1\lambda^p_k})^F$
coincides with the essential spectrum of the Friedrichs extension of the
Laplacian $-\frac{\partial^2}{\partial r^2}$ on
$C^{\infty}_c(\tilde c,+\infty)$.
Hence, if $\lambda^p_k=0$, $\sigma_{\rm ess}((D_{1\lambda^p_k})^F)=\left[
0,+\infty
\right)$.\par
If, on the contrary, $\lambda^p_k >0$, for every $h \in (0,1)$,
$$ \lim_{r \rightarrow +\infty}\int_r^{r+h}
\left(K_1(n,p)s^{-2}+\lambda^p_k
|a+1|^{\frac{2b}{|a+1}}s^{-\frac{2b}{a+1}} \right)\,ds=$$
$$ = \lim_{r\rightarrow +\infty} \lambda^p_k \left(\frac{4b}{|a+1|}+1
\right)^{-1}\left[ (r+h)^{\frac{2b}{|a+1|}+1} -r^{\frac{2b}{|a+1|}+1}
\right]=+\infty$$
since $\frac{2b}{|a+1|}>0$. Then, in view of
Theorem 3.13 in \cite{Muller-Pfeiffer},
$\sigma_{\rm ess}(D_{1\lambda^p_k}) = \emptyset$.
\end{proof}
If $K_1(n,p)\geq 0$, we find
$$ \langle (D_{1\lambda^p_k})^F w,w\rangle_{L^2(\tilde c,+\infty)} \geq
\int_{\tilde c}^{+\infty} \lambda^p_k |a+1|^{\frac{2b}{|a+1|}}
s^{\frac{2b}{|a+1|}}
w^2 \,ds \geq$$
$$ \geq  C |a+1|^{\frac{2b}{|a+1|}} \lambda^p_k \langle
w,w\rangle_{L^2(\tilde c,+\infty)}.$$
If $K_1(n,p)<0$, since, again, $\tilde c>0$ can be chosen
arbitrarily large,
we can suppose that for every $r>\tilde c$ $$V(r):= K_1(n,p) r^{-2}+
\lambda^p_k
|a+1|^{\frac{2b}{|a+1}}r^{\frac{2b}{|a+1|}}> V(\tilde c)>0;$$
hence,
\begin{multline}\label{stima<-1>0} \langle (D_{1\lambda^p_k})^F w, w
\rangle_{L^2(\tilde c,+\infty)} \geq V(\tilde c)
\langle w,w
\rangle_{L^2(\tilde c,+\infty)}=\\
= \left( \frac{K_1(n,p)}{{\tilde c}^2} + \lambda^p_k
|a+1|^{\frac{2b}{|a+1|}}
{\tilde  c}^{\frac{2b}{|a+1|}} \right)\langle
w,w\rangle_{L^2(\tilde c,+\infty)}.\end{multline}
In both cases, for any $w$ such that $\|w\|_{L^2(\tilde c,+\infty)}=1$
$$\langle (D_{1\lambda^p_k})^Fw,w
\rangle_{L^2(\tilde c,+\infty)}\rightarrow +\infty $$ as $k\rightarrow
+\infty$. As a consequence:
\begin{prop}\label{<-1>01M} Let $a<-1$, $b>0$. For $0\leq p \leq
(n-1)$, if the
$p$-th Betti number of $N$ vanishes, then $\sigma_{\rm
ess}((\Delta^p_{M1})^F)= \emptyset$. If, on the contrary, the $p$-th
Betti number of $N$ is different from zero, then $\sigma_{\rm
ess}((\Delta^p_{M1})^F)= \left[ 0,+\infty\right)$.
\end{prop}
By duality,
\begin{prop}\label{<-1>02M} Let $a<-1$, $b>0$. For $1\leq p \leq
n$, if the $(p-1)$-th Betti number of $N$ vanishes, then $\sigma_{\rm
ess}((\Delta^p_{M2})^F)= \emptyset$; if, on the contrary, the
$(p-1)$-th Betti number of $N$ is different from zero, then $\sigma_{\rm
ess}((\Delta^p_{M2})^F)= \left[ 0, +\infty\right)$.
\end{prop}
We still have to investigate the essential spectrum of
$(\Delta^p_{M3})^F$. First of all, we will compute the essential
spectrum of $D_{3\lambda_k^{p-1}}$ for every $k\in \mathbb N$. In
analogy with what we did in the case $a=-1$, $b>0$, we need a
preliminary Lemma:
\begin{lem}\label{fittizio<-1} For every $K \in \Real$, the essential
spectrum of the
Friedrichs extension $D^F$ of the operator
$$D:C^{\infty}_c(c,+\infty)\oplus C^{\infty}_c(c,+\infty) \longrightarrow
L^2(c,+\infty)\oplus L^2(c,+\infty) $$
defined by
$$D(w_1 \oplus w_2)= \left(- \frac{\partial^2 w_1}{\partial
r^2}+Kr^{-2}w_1 +
\lambda^{p-1}_k|a+1|^{-\frac{2b}{a+1}}\,r^{-\frac{2b}{a+1}} w_1+ \right.$$
$$\left. + \sqrt{\lambda^{p-1}_k}\,|a+1|^{-\frac{b}{a+1}}\,
r^{-\frac{b}{a+1}-1} w_2 \right) \oplus$$
$$ \oplus \left( - \frac{\partial^2 w_2}{\partial r^2}+ K r^{-2}w_2 +
\lambda^{p-1}_k|a+1|^{-\frac{2b}{a+1}}\,r^{-\frac{2b}{a+1}}w_2 +\right.$$
$$\left. +
\sqrt{\lambda^{p-1}_k}|a+1|^{-\frac{b}{a+1}}\,r^{-\frac{b}{a+1}-1}\right)$$
is empty.
\end{lem}
\begin{proof} First of all, since as previously stated the essential spectrum of $D^F$ does not
depend on the choice of the first endpoint $\tilde c$ of $(\tilde c,
+\infty)$, given $K$ and $\mu >0$ we can suppose that
for some positive constant $C$, for every $t>\tilde c$
\begin{multline}\label{estest}K r^{-2}+
\lambda^{p-1}_k|a+1|^{2\frac{|b|}{|a+1|}}r^{2\frac{|b|}{|a+1|}}
\pm
\sqrt{\lambda^{p-1}_k}|a+1|^{\frac{|b|}{|a+1|}}r^{\frac{|b|}{|a+1|}-1}-\mu
\geq C.\end{multline}
As in the proof of Lemma \ref{fittizio}, let us consider the orthogonal
decomposition (\ref{decc})
$$L^2(\tilde c,+\infty)\oplus L^2(\tilde c,+\infty)= {\mathcal V}_1\oplus
{\mathcal
V}_2,
$$
where
$${\mathcal V}_1:= \left\{w_1\oplus w_2 \,|\, w_1=w_2 \right\} $$
and
$${\mathcal V}_2:= \left\{ w_1\oplus w_2 \,|\, w_1=-w_2 \right\}.$$
$D$ is invariant under (\ref{decc}), hence
$$ D^F= (D_{|{\mathcal V}_1})^F\oplus(D_{|{\mathcal V}_2})^F,$$
and
$$\sigma_{\rm ess}(D^F)= \sigma_{\rm ess}(D_{|{\mathcal V}_1}^F)\cup
\sigma_{\rm ess}(D_{|{\mathcal V}_2}^F).$$
If $\mu \in \sigma_{\rm ess}(D_{|{\mathcal V}_1}^F)$, there exists a Weyl
sequence $\left\{w_n\right\}\subset {\mathcal D}(D_{|{\mathcal V}_1^F})$
for $\mu$: we have
that $\|w_n\|_{L^2(\tilde c,+\infty)}\leq C$, $(D-\mu)(w_n\oplus
w_n)\longrightarrow 0 $ as $n \rightarrow +\infty$, but
$\left\{w_n\right\}$ has no convergent subsequence. Then
$$\langle (D -\mu) (w_n\oplus w_n),(w_n\oplus w_n)
\rangle_{L^2(\tilde c,+\infty)\oplus
L^2(\tilde c,+\infty)}\longrightarrow 0 $$
as $n\rightarrow +\infty$; hence
\begin{multline}\int_{\tilde c}^{+\infty}\left(\frac{\partial
w_n}{\partial
x}\right)^2\,dx+
\left( K x^{-2}+
\lambda^{p-1}_k|a+1|^{2\frac{|b|}{|a+1|}}x^{2\frac{|b|}{|a+1|}}+
\right.\\
 -
\left.
\sqrt{\lambda^{p-1}_k}|a+1|^{\frac{|b|}{|a+1|}}x^{\frac{|b|}{|a+1|}-1}
-\mu\right)w_n^2\,dx \longrightarrow 0\end{multline} as
$n\rightarrow +\infty$. Since the estimate (\ref{estest}) implies
that $\|w_n\|_{L^2(\tilde c,+\infty)} \rightarrow 0$ as
$n\rightarrow +\infty$, we obtain that $\mu \notin \sigma_{\rm
ess}(D^F_{|{\mathcal V}_1})$. \par Suppose now that $\mu \in
\sigma_{\rm ess}(D^F_{|{\mathcal V}_2})$. If $\left\{w_n\right\}$
is a Weyl sequence for $\mu$ for $D_{|{\mathcal V}_2}^F$, we have
$$\langle (D-\mu)(w_n\oplus -w_n),(w_n\oplus -
w_n)\rangle_{L^2(\tilde c,+\infty)\oplus L^2(\tilde c,+\infty)}
\longrightarrow 0$$ as $n\rightarrow +\infty$. Again in view of
the estimate (\ref{estest}), $\|w_n\|_{L^2(\tilde
c,+\infty)}\rightarrow 0$ as $n \rightarrow +\infty$. But a Weyl
sequence cannot converge. Hence, $\mu \notin \sigma_{\rm
ess}(D^F_{|{\mathcal V}_2})$. This completes the proof.\end{proof}
We can now compute the essential spectrum of
$(D_{3\lambda_k^{p-1}})^F$ for every $k\in \mathbb N$:
\begin{lem} For $1\leq p \leq n-1$ and for every $k \in \mathbb N $,
$$\sigma_{\rm ess}((D_{3\lambda_k^{p-1}})^F)= \emptyset.$$
\end{lem}
\begin{proof} Let us consider the Friedrichs extension $(D_{30})^F$ of the
operator $$D_{30}:C^{\infty}_c(\tilde c,+\infty)\oplus
C^{\infty}_c(\tilde c,+\infty)\longrightarrow L^2(\tilde
c,+\infty)\oplus L^2(\tilde c,+\infty) $$ $$D_{30}(w_1\oplus w_2)=
\left( - \frac{\partial^2 w_1}{\partial r^2}+ \bar K r^{-2}w_1
+\lambda^{p-1}_k|a+1|^{-\frac{2b}{a+1}}r^{-\frac{2b}{a+1}}w_1 +
\right. $$ $$\left.
+\sqrt{\lambda^{p-1}_k}|a+1|^{-\frac{b}{a+1}}\,r^{-\frac{b}{a+1}-1}w_2
\right) \oplus $$ $$\oplus \left( - \frac{\partial^2 w_2}{\partial
r^2}+ \bar K r^{-2}w_2
+\lambda^{p-1}_k|a+1|^{-\frac{2b}{a+1}}\,r^{-\frac{2b}{a+1}}w_2
+\right.$$ $$ \left.
+\sqrt{\lambda^{p-1}_k}|a+1|^{-\frac{b}{a+1}}\,r^{-\frac{b}{a+1}-1}w_1\right),
$$ where $$ \bar K= \max \left\{K_1(n,p),K_2(n,p) \right\}.$$ From
the previous Lemma, we know that $\sigma_{\rm ess}((D_{30})^F)=
\emptyset$. We will show that
$((D_{3\lambda_k^{p-1}})^F-(D_{30})^F)$ is a relatively compact
perturbation of $(D_{30})^F$. First of all, a straightforward
computation shows that $X_{D_{30}}\subseteq
X_{D_{3\lambda_k^{p-1}}}$ and ${\mathcal D}(D_{30}^*)\subseteq
{\mathcal D}(D_{3\lambda_k^{p-1}}^*)$; hence, ${\mathcal
D}((D_{30})^F)\subseteq {\mathcal D}((D_{3\lambda_k^{p-1}}^F))$,
whence ${\mathcal D}((D_{3\lambda_k^{p-1}})^F-(D_{30})^F)=
{\mathcal D}((D_{30})^F)$.\par We still have to show that, given a
sequence $\left\{ w_{1n}\oplus w_{2n}\right\}$ \\ $\subset
{\mathcal D}((D_{30})^F)$ such that $$\|w_{1n}\oplus
w_{2n}\|^2_{L^2\oplus L^2}+\|((D_{30})^F)(w_{1n}\oplus
w_{2n})\|^2_{L^2\oplus L^2}\leq C,$$ there exists a subsequence
$\left\{ w_{1n_k}\oplus w_{2n_k}\right\}$ such that
$$\left\{((D_{3\lambda_k^{p-1}})^F-(D_{30})^F)(w_{1n_k}\oplus
w_{2n_k})\right\}$$ converges.\par Now, the fact that
$\|(D_{30})^F(w_{1n}\oplus w_{2n})\|_{L^2(\tilde c, +\infty)\oplus
L^2(\tilde c, +\infty)}\leq C$ is equivalent to the inequalities
\begin{multline}\label{<uno} \|- \frac{\partial^2 w_{1n}}{\partial r^2}+
\bar
K r^{-2}w_{1n}
+\lambda^{p-1}_k|a+1|^{-\frac{2b}{a+1}}r^{-\frac{2b}{a+1}}w_{1n}+
\\+\sqrt{\lambda^{p-1}_k}|a+1|^{-\frac{b}{a+1}}\,r^{-\frac{b}{a+1}-1}
w_{2n}\|_{L^2(\tilde c,+\infty)}\leq
C,
\end{multline}
\begin{multline}\label{<due}\| - \frac{\partial^2 w_{2n}}{\partial r^2}+
\bar K
r^{-2}w_{2n}
+\lambda^{p-1}_k|a+1|^{-\frac{2b}{a+1}}r^{-\frac{2b}{a+1}}w_{2n}+
\\+\sqrt{\lambda^{p-1}_k}|a+1|^{-\frac{b}{a+1}}\,r^{-\frac{b}{a+1}-1}w_{1n}
\|_{L^2(\tilde c,+\infty)}\leq C,
\end{multline}
which in turn imply
\begin{multline}\label{<+} \|- \frac{\partial^2 (w_{1n}+w_{2n})}{\partial
r^2}+ \bar K r^{-2} (w_{1n}+w_{2n})+\\
+\lambda^{p-1}_k|a+1|^{-\frac{2b}{a+1}}r^{-\frac{2b}{a+1}}(w_{1n}+w_{2n})+\\
+ \sqrt{\lambda^{p-1}_k} |a+1|^{-\frac{b}{a+1}}\,r^{-
\frac{b}{a+1}-1}(w_{1n}+w_{2n})\|_{L^2(\tilde c, +\infty)}\leq C,
\end{multline} \begin{multline}\label{<-}\|- \frac{\partial^2
(w_{1n}-w_{2n})}{\partial
r^2}+ \bar K r^{-2} (w_{1n}-w_{2n})+\\
+\lambda^{p-1}_k|a+1|^{-\frac{2b}{a+1}}r^{-\frac{2b}{a+1}}(w_{1n}-w_{2n})+\\
+ \sqrt{\lambda^{p-1}_k} |a+1|^{-\frac{b}{a+1}}\,r^{-
\frac{b}{a+1}-1}(w_{2n}-w_{1n})\|_{L^2(\tilde c, +\infty )}\leq C.
\end{multline}
By multiplication with $(w_{1n}+w_{2n})$ and $(w_{1n}-w_{2n})$
respectively, (\ref{<+}) and (\ref{<-}) yield:
\begin{multline}\label{<++}
\int_{\tilde c}^{+\infty}\left(\frac{\partial(w_{1n}+w_{2n})}{\partial
s}\right)^2\,ds + \int_{\tilde c}^{+\infty}  \bar K
s^{-2}
(w_{1n}+w_{2n})^2\,ds\\
+\int_{\tilde c}^{+\infty}\lambda^{p-1}_k|a+1|^{-\frac{2b}{a+1}}
s^{-\frac{2b}{a+1}}
(w_{1n}+w_{2n})^2\,ds+\\
+ \int_{ \tilde c}^{+\infty} \sqrt{\lambda^{p-1}_k}
|a+1|^{-\frac{b}{a+1}}\,s^{-
\frac{b}{a+1}-1}(w_{2n}+w_{1n})^2\,ds
\leq C, \end{multline}
\begin{multline}\label{<--}\int_{\tilde c}^{+\infty}\left(
\frac{\partial(w_{1n}-w_{2n})}{\partial
s}\right)^2\,ds + \int_{\tilde c}^{+\infty} \bar K
s^{-2}
(w_{1n}-w_{2n})^2\,ds +
\\
+\int_{\tilde c}^{+\infty}\lambda^{p-1}_k|a+1|^{-\frac{2b}{a+1}}
s^{-\frac{2b}{a+1}}
(w_{1n}-w_{2n})^2\,ds+\\
+ \int_{/tilde c}^{+\infty} - \sqrt{\lambda^{p-1}_k}
|a+1|^{-\frac{b}{a+1}}\,s^{-
\frac{b}{a+1}-1}(w_{1n}-w_{2n})^2\,ds \leq C.
\end{multline}
Now, there exists $\bar c$ such that for every $k \in \mathbb N$
and for every $r>\bar c$ $$\bar K r^{-2}
+\lambda^{p-1}_k|a+1|^{\frac{|2b|}{|a+1|}}r^{\frac{|2b|}{|a+1|}}\pm
\sqrt{\lambda^{p-1}_k}|a+1|^{\frac{|2b|}{|a+1|}}r^{\frac{|b|}{|a+1|}-1}
\geq C $$ for some positive constant $C$. Since, in view of the
results of \cite{Eichhorn}, the essential spectrum of
$(D_{3\lambda_k^{p-1}})^F$ does not depend on the first endpoint
$\tilde c$ of $(\tilde c, +\infty)$, we can suppose that $\tilde
c>\bar c$; then,  from the estimates (\ref{<++}) and (\ref{<--})
we obtain that $$\int_{\tilde c}^{+\infty}  (w_{1n}\pm
w_{2n})^2\,ds \leq C. $$ Hence, $\left\{w_{1n}\right\}$ and
$\left\{w_{2n}\right\}$ are bounded in $W^{1,2}(\tilde
c,+\infty)$; as a consequence, they are bounded also in
$L^{\infty}(\tilde c,+\infty)$ and in $W^{1,2}(K)$ for every
compact subset $K \subset (\tilde c,+\infty)$. Now, for every
$n,m\in \mathbb N$ $$
\|((D_{3\lambda^{p-1}_k})^F-(D_{30})^F)((w_{1n}-w_{1m})\oplus(w_{2n}-w_{2m}))
\|_{L^2(\tilde c, +\infty)\oplus L^2(\tilde c, +\infty)}^2 =$$ $$=
\sum_{i=1,2}\int_{\tilde c}^{+\infty}(K_i(p,n)-\bar K)^2 s^{-4}
(w_{in}-w_{im})^2\,ds;$$ since $s^{-2}\in L^2(\tilde
c,+\infty)\cap L^{\infty}(\tilde c,+\infty)$, following the
argument of Lemma \ref{-1<01l} the conclusion follows.
\end{proof}
Moreover,
\begin{lem}\label{ficicchiosa} Let $a<-1$, $b>0$. If the $p$-th and the
$(p-1)$-th Betti numbers of $N$ both
vanish, and if $0<\mu \in \overline{\sigma_p((\Delta^p_{M3})^F)}$, then
$\mu$ is an isolated eigenvalue of finite multiplicity.
\end{lem}
\begin{proof} If $\mu$ is an eigenvalue of infinite multiplicity of
$(\Delta^p_{M3})^F$ or is a cluster point of
$\sigma_p((\Delta^p_{M3})^F)$, there
exist a sequence $\left\{\mu_k \right\}$ of eigenvalues of
$(\Delta^p_{M3})^F$ and
a corresponding sequence of normalized, mutually orthogonal
eigenforms $\left\{\Phi_k \right\}$ such that for every $k$
$$ (\Delta^p_{M3})^F \Phi_k - \mu_k \Phi_k =0$$
and
$$ \mu_k \longrightarrow \mu \quad \quad \hbox{as $k \rightarrow
+\infty$.}$$
In view of the weak Kodaira decomposition, replacing
$\left\{\Phi_k\right\}$ by a subsequence (again denoted, for shortness, by
the same symbol) we can suppose that either $\delta^p_M \Phi_k =0$
for every $k\in \mathbb N$ or $d^p_M \Phi_k =0$ for every $k\in \mathbb
N$. Following the argument of Lemma
\ref{difficile}, we find that in the first case $\mu$ lies in the
essential spectrum of the operator $(\Delta^{p+1}_{M2})^F$,
 whilst in the second case $\mu$ lies in the essential
spectrum of the operator $(\Delta^{p-1}_{M1})^F$. \par \noindent
In the first case, since by assumption the $p$-th Betti number
of $N$ vanishes, we get a contradiction with Proposition \ref{<-1>02M}.
In the second case, since by assumption the $(p-1)$-th Betti number
of $N$ vanishes, we get a contradiction with Proposition \ref{<-1>01M}.
\end{proof}
Summing up, for $a<-1$, $b>0$ we have the following result:
\begin{thm} Let $a<-1$, $b>0$, $p\in [0,n]$. If both the $p$-th and the
$(p-1)$-th Betti
numbers of $N$ vanish, then $\sigma_{\rm ess}(\Delta^p_M)\setminus
\left\{0\right\}= \emptyset$. If, on the contrary, at least one of them is
different from zero,
then $ \sigma_{\rm ess}(\Delta^p_M)=\left[0,+\infty\right)$.
\end{thm}
\section{The rotationally symmetric case: $N= \Sphere^{n-1}$}
Since the techniques in Lemma \ref{difficile}, \ref{ficicchia},
\ref{ficicchiosa} work only when $\mu_k >0$, in the case of a
general manifold with boundary $\overline M$ for some values of
the parameters $a$ and $b$ (namely for $a=-1$, $b<0$, for $a=-1$,
$b>0$ and for $a<-1$, $b>0$) we were able to compute only
$\sigma_{\rm ess}(\Delta^p_M)\setminus \left\{ 0 \right\}$. In
particular, in those cases we were not able to establish whether
$0$ is an eigenvalue of infinite multiplicity of $\Delta_M^p$ or
not. \par In the rotationally symmetric case, that is, when
$\overline M $ is the unitary ball $\overline{ B(\underline 0,1)}$
in $\Real^n$ and the Riemannian metric is globally invariant under
rotations, we can be slightly more precise. In fact, we have the
following generalization of a result of Dodziuk (see
\cite{Dodziuk}), whose proof essentially follows that of
\cite{Dodziuk} and is therefore omitted.
\begin{thm}\label{Dodziuk1}
Let us consider, for $n\geq 2$, the manifold $$\tilde M=[0,+\infty)\times
\Sphere^{n-1}$$ endowed
with a complete Riemannian metric of the type
\begin{equation}\label{type}f(t)\,dt^2+g(t)\,d\theta^2_{\Sphere^{n-1}},
\end{equation}
where $f(t)>0$, $g(t)>0$ for every $t\in [0,+\infty)$, and
$d\theta_{\Sphere^{n-1}}^2$ is the standard metric on $\Sphere^{n-1}$.
Let us suppose that for some $\epsilon >0$
\begin{equation}\label{condzero} f(t)\equiv 1,\quad g(t)=t^2 \quad
\hbox{for $t\in (0,\epsilon)$}.\end{equation}
If
we denote by ${\mathcal H}^p(\tilde M)$, for $p=0,...,n$,
the space of $L^2$ harmonic $p$-forms on $\tilde M$, we have:
\begin{enumerate}
\item for $p\not=0,n,n/2$, ${\mathcal H}^p(\tilde M)=\left\{0\right\}$;
\item if $\int_0^{\infty}f^{\frac{1}{2}}(s)g^{\frac{n-1}{2}}(s)\,ds=
+\infty$, then ${\mathcal H}^n( \tilde M)\simeq {\mathcal H}^0(\tilde
M)=\left\{0\right\}$; if, on the contrary,
$\int_0^{\infty}f^{\frac{1}{2}}(s)g^{\frac{n-1}{2}}(s)\,ds<+
\infty$, then
${\mathcal H}^n(\tilde M)\simeq {\mathcal H}^0(\tilde M)=\Real$;
\item for $p=\frac{n}{2}$, if $\int_1^{+\infty}
f^{\frac{1}{2}}(s)g^{-\frac{1}{2}}(s)\,ds = +\infty$,
then ${\mathcal H}^p(\tilde M)= \left\{0\right\}$; if,
on the other hand,
$\int_1^{+\infty}f^{\frac{1}{2}}(s)g^{-\frac{1}{2}}(s)\,ds < +\infty$,
then ${\mathcal H}^{\frac{n}{2}}(\tilde M)$ is a Hilbert space of infinite
dimension.
\end{enumerate}
\end{thm}
Now, suppose that the metric (\ref{type}) fulfills the asymptotic
condition:
\begin{equation}\label{condinf}f(t)= e^{-2(1+a)t}, \quad g(t)= e^{-2bt}
\quad \hbox{for $t>c>>0$}. \end{equation}
A straightforward computation shows that
$$\int_0^{+\infty}f(s)^{\frac{1}{2}}g(s)^{\frac{n-1}{2}}\,ds \left\{
\begin{array}{ll} = +\infty \quad &\hbox{if $b\leq - \frac{a+1}{n-1}$}\\
< +\infty \quad &\hbox{if $b> - \frac{a+1}{n-1}$}
\end{array}, \right. $$
whilst
$$\int_1^{+\infty}f(s)^{\frac{1}{2}}g(s)^{-\frac{1}{2}}\,ds \left\{
\begin{array}{ll}=+\infty \quad & \hbox{if $b\geq a+1$}\\
<+\infty \quad & \hbox{if $b<a+1$}
\end{array}.\right. $$
As a consequence, from Theorem \ref{Dodziuk1} we can easily deduce the
following
\begin{thm}\label{Dodziuk2} For $n\geq 2$,  let us consider the manifold
$\tilde M$,
endowed with
a Riemannian metric of type (\ref{type}), satisfying conditions
(\ref{condzero}) and (\ref{condinf}).
\begin{enumerate}
\item If $p\not=n/2,n,0$, then $0\notin
\sigma_p(\Delta^p_{\tilde M})$.
\item If $p=0,n$, then $0\in \sigma_p(\Delta^p_{\tilde M})$ if and
only if $b> -
\frac{a+1}{n-1}$.
\item If $p=n/2$, $0\in \sigma_p(\Delta^p_{\tilde M})$ if and only if
$b<a+1$, and in
this case ${\mathcal H}^p(\tilde M)$ is a Hilbert space of infinite
dimension, hence $0\in \sigma_{{\rm ess}}(\Delta^p_{\tilde M})\cap
\sigma_p(\Delta^p_{\tilde M})$.
\end{enumerate}
\end{thm}
Let us consider, for the moment, the Riemannian manifold $\tilde
M$. Since the Riemannian metric on $\tilde M$ is of type
(\ref{fg}), the results of sections 4 and 5 can be recovered
through the decomposition techniques of section 3. Moreover, since
the Riemannian metric is globally invariant under rotations, the
decompositions can be applied directly on $\Delta^p_{\tilde M}$
and not only on its Friedrichs extension $(\Delta^p_{\tilde
M})^F$: we can write $$ \Delta^p_{\tilde M}=\Delta^p_{\tilde
M1}\oplus \Delta^p_{\tilde M2} \oplus \Delta^p_{\tilde M3},$$
where, for $i=1,2,3$, $\Delta^p_{\tilde M i}$ is the restriction
of $\Delta_{\tilde M}^p$ to ${\mathcal L}_i(M)$". Moreover, we
have $$\sigma_{\rm ess}(\Delta^p_M)=\bigcup_{i=1,2,3}\sigma_{\rm
ess}(\Delta^p_{Mi}) ,$$ $$\sigma_p(\Delta_p^M)=
\sigma_p(\Delta^p_{Mi}) .$$ In the proof of Proposition
\ref{difficile} we used the fact that if $\left\{\Phi_k\right\}$
is a sequence of normalized, mutually orthogonal $p$-eigenforms of
$(\Delta^p_{M3})^F$ corresponding to some positive eigenvalues
$\mu_k$ and $\mu_k \rightarrow \mu>0$ as $k\rightarrow +\infty$,
then either $\mu$ is in the essential spectrum of
$(\Delta^{p-1}_{M1})^F$ or $\mu$ is in the essential spectrum of
$(\Delta^{p+1}_{M2})^F$. In the rotationally symmetric case, this
property can be considerably strengthened; namely, the following
Proposition holds:
\begin{prop}\label{link} If $\Phi$ is a $p$-eigenform for
$\Delta^p_{\tilde M3}$
corresponding to a positive eigenvalue $\mu$, if $d^p_{\tilde
M}\Phi\not=0$
(resp. $\delta^p_{\tilde M}\Phi \not= 0$), then $d^p_{\tilde M}\Phi$
(resp.
$\delta^p_{\tilde M}\Phi$) is a $(p+1)$-eigenform of $\Delta^{p+1}_{\tilde
M2}$
(resp. a $(p-1)$-eigenform of $\Delta^{p-1}_{\tilde M1}$).
 \end{prop}
Combining the results of sections 4, 5 with Theorem \ref{Dodziuk2} and
with Proposition \ref{link} we find:
\begin{prop}\label{zero} Let $\tilde M$ be as in Theorem \ref{Dodziuk2}:
then
\begin{enumerate}
\item if $a=-1$, $b<0$, $0 \in \sigma_{\rm ess}(\Delta^p_{\tilde M})$ if
and
only if $p \in \left\{ \frac{n}{2},\frac{n+ 1}{2}, \frac{n-1}{2}
\right\}$;
\item if $a=-1$, $b>0$, for every $p \in [0,n]$ $0 \notin \sigma_{\rm
ess}(\Delta^p_{\tilde M})$;
\item if $a<-1$ and $b>0$, $0 \in \sigma_{\rm ess}(\Delta^p_{\tilde M})$
if
and only if \\ $p \in \left\{ 0,1,n-1,n \right\}$.
\end{enumerate}
\end{prop}
\begin{proof} 1. Let $a=-1$, $b<0$. For $p=\frac{n\pm 1}{2}$, $0\in
\sigma_{\rm ess}(\Delta^p_{\tilde M})$ thanks to Theorem
\ref{a=-1b<0}. On the other hand, if $p=\frac{n}{2}$, $0$ is an
eigenvalue of infinite multiplicity by Theorem \ref{Dodziuk2},
hence $0\in \sigma_{\rm ess}(\Delta^p_{\tilde M})$. If, on the
contrary, $p\notin \left\{ \frac{n-1}{2}, \frac{n}{2} ,
\frac{n+1}{2}\right\}$, then $0\notin \sigma_p(\Delta^p_{\tilde
M})$ by Theorem \ref{Dodziuk2}; hence, in view of Proposition
\ref{-1<0M1}, if $0\in \sigma_{\rm ess}(\Delta^p_{\tilde M})$
there exists a sequence $\left\{\mu_k\right\}$ of positive
eigenvalues of $\Delta^p_{\tilde M3}$ and a corresponding sequence
of normalized, mutually orthogonal $p$-eigenforms
$\left\{\Phi_k\right\}$ of $\Delta^p_{M3}$ such that $\mu_k
\rightarrow 0$ as $k\rightarrow +\infty$. Thanks to Proposition
\ref{link}, either $0\in \sigma_{\rm ess}(\Delta^{p-1}_{\tilde
M1})$ or $0\in \sigma_{\rm ess}(\Delta^{p+1}_{\tilde M2})$. Then,
in view of Proposition \ref{-1<0M2}, we should have $$ 0\geq
\left(\frac{n-2(p-1)-1}{2}\right)^2b^2=
\left(\frac{n-2p+1}{2}\right)^2b^2$$ or $$0\geq
\left(\frac{n-2(p+1)+1}{2} \right)^2 b^2
=\left(\frac{n-2p-1}{2}\right)^2b^2;$$ since $p\not= \frac{n\pm
1}{2}$, we get a contradiction.\par \noindent

2. Let $a=-1$, $b>0$. Since $N= \Sphere^{n-1}$, the $p$-th Betti
number of $N$ is different from zero if and only if $p=0$ or
$p=n-1$. As a consequence, in view of Propositions \ref{betti1}
and \ref{betti2}, we obtain that, for
$p\notin\left\{0,1,n-1,n\right\}$, $\sigma_{\rm
ess}(\Delta^p_{\tilde M1})= \sigma_{\rm ess}(\Delta^p_{\tilde
M2})= \emptyset$. Since Theorem \ref{Dodziuk2} implies that for
$p\notin \left\{0,1,n-1,n\right\}$ $0$ cannot be an eigenvalue of
$\Delta_{\tilde M}^p$, $0$ can lie in the essential spectrum of
$\Delta^p_{\tilde M}$ if and only if there exists a sequence
$\left\{\mu_k\right\}$ of positive eigenvalues of
$\Delta^p_{\tilde M3}$ and a corresponding sequence $\left\{
\Phi_k \right\}$ of normalized, mutually orthogonal $p$-eigenforms
of $\Delta^p_{\tilde M3}$ such that $\mu_k\rightarrow 0$ as
$k\rightarrow +\infty$. Then, Proposition \ref{link} implies that
either $0\in\sigma_{\rm ess}(\Delta^{p-1}_{\tilde M1})$ or $0\in
\sigma_{\rm ess}(\Delta^{p+1}_{\tilde M2})$. Since $p-1 \notin
\left\{0,n-1\right\}$ and $(p+1)-1\notin \left\{0,n-1 \right\}$,
we get a contradiction with Proposition \ref{betti1} or with
Proposition \ref{betti2}. Hence, if $p\notin \left\{0,1,n-1,n
\right\}$, then $0\notin \sigma_{\rm ess}(\Delta^p_{\tilde M})$.
\par
If $p=0$, we have only $p$-forms of type I and, since $n>1$,
$0\notin \sigma_{\rm ess}(\Delta^p_{\tilde M})$. By duality, the
same result holds for $p=n$.\par Finally, let $p=1$; Propositions
\ref{betti1} and \ref{betti2} imply that $\sigma_{\rm
ess}(\Delta^p_{\tilde M1})= \emptyset$ whilst $$\sigma_{\rm
ess}(\Delta^p_{\tilde M2})= \left[\left(\frac{n-1}{2}\right)^2b^2
,+\infty \right).$$ Since, in view of Theorem \ref{Dodziuk2}, $0$
cannot be an eigenvalue, following the same argument as before we
see that if $0\in \sigma_{\rm ess}(\Delta^1_{\tilde M})$, then
either $0\in \sigma_{\rm ess}(\Delta^0_{\tilde M 1})$ or $0\in
\sigma_{\rm ess}(\Delta^2_{\tilde M 2})$. But this is in
contradiction with Proposition \ref{betti1} or with Proposition
\ref{betti2}. Hence, $0\notin \sigma_{\rm ess}(\Delta^p_{\tilde
M})$ for any $p\in [0,n]$.\par \noindent

3. Let $a<-1$ and $b>0$; again, since $N= \Sphere^{n-1}$, the $p$-th
Betti number of $N$ is different from zero if and only if $p=0$ or
$p=n-1$. \par
Then, if $p\notin \left\{0,1,n-1,n \right\}$, by Propositions
\ref{<-1>01M}
and \ref{<-1>02M} $$\sigma_{\rm ess}(\Delta^p_{\tilde M1})= \sigma_{\rm
ess}(\Delta^p_{\tilde M2})= \emptyset.$$ Since, by Theorem \ref{Dodziuk2},
$0$
cannot be an eigenvalue of $\Delta^p_{\tilde M}$, if $0\in \sigma_{\rm
ess}(\Delta^p_{\tilde M})$ there exist a sequence $\left\{\mu_k\right\}$
of
positive eigenvalues of $\Delta^p_{\tilde M3}$ and a corresponding
sequence
$\left\{ \Phi_k \right\}$ of normalized, mutually orthogonal eigenforms of
$\Delta^p_{\tilde M3}$ such that $\mu_k \rightarrow 0$
as $k\rightarrow +\infty$.
Then, by Proposition \ref{link}, either $0\in \sigma_{\rm
ess}(\Delta^{p-1}_{\tilde M1})$  or $0\in \sigma_{\rm
ess}(\Delta^{p+1}_{\tilde M2})$. Since $1<p<n-1$, we get
a
contradiction with Propositions \ref{<-1>01M} and \ref{<-1>02M}. \par
If, on the contrary, $p\in \left\{ 0,1,n-1,n \right\}$, again by
Propositions \ref{<-1>01M} and \ref{<-1>02M}, $0\in \sigma_{\rm
ess}(\Delta_{\tilde M})$.
\end{proof}
Since the essential spectra of $\tilde M$ and of $M$ coincide,
combining the results of section 4 with  Proposition \ref{zero},
we finally get:
\begin{thm} Let $M$ be the unitary ball $B(\underline 0,1)$ in $\Real^n$
endowed with a Riemannian metric which, in a tubular neighbourhood of the
boundary $\Sphere^{n-1}$, is given by
$$d\sigma^2= e^{-2(a+1)t}\,dt^2 + e^{-2bt}\,d\theta^2_{\Sphere^{n-1}},$$
where $a\leq -1$, $t=settanh(\|x\|) $ and $d\theta^2_{\Sphere^{n-1}}$ is
the standard
Riemannian
metric on $\Sphere^{n-1}$. Then
\begin{enumerate}
\item if $a=-1$ and $b<0$, if $p \not= \frac{n}{2}$ $$\sigma_{\rm
ess}(\Delta^p_M)= \left[\min\left\{
\left(\frac{n-2p-1}{2}\right)^2b^2,\left(\frac{n-2p+1}{2}\right)^2b^2\right\}
,+\infty\right), $$
whilst if $p=\frac{n}{2}$
$$\sigma_{\rm ess}(\Delta^p_M)= \left\{0\right\}\cup \left[ \frac{b^2}{4},
+\infty \right) ;$$
\item if $a=-1$ and $b=0$, then for every $p\in [0,n]$ $$\sigma_{\rm
ess}(\Delta^p_M)=[\overline \lambda,+\infty),$$
where $\overline \lambda$ is the minimum between the smallest eigenvalue
of
$\Delta^p_{\Sphere^{n-1}}$ and the smallest eigenvalue
of $\Delta^{p-1}_{\Sphere^{n-1}}$;
\item if $a=-1$ and $b>0$, if $1<p<n-1$ $$\sigma_{\rm ess}(\Delta^p_M)=
\emptyset,$$ whilst if $p\in \left\{0,1,n-1,n \right\}$ $$\sigma_{\rm
ess}(\Delta^p_M)= \left[ \left(\frac{n-1}{2}\right)^2 b^2,
+\infty\right);$$
\item if $a<-1$ and $b<0$, then, for every $p\in [0,n]$, $\sigma_{\rm
ess}(\Delta^p_M)=[0,+\infty);$
\item if $a<-1$ and $b=0$, for every $p\in [0,n]$ $$\sigma_{\rm
ess}(\Delta^p_M)=[\overline\lambda,+\infty),$$where $\overline\lambda$
is the minimum
between the smallest eigenvalue of
$\Delta^p_{\Sphere^{n-1}}$ and the smallest eigenvalue
of $\Delta^{p-1}_{\Sphere^{n-1}}$;
\item if $a<-1$ and $b>0$, if $1<p<n-1$ then $$\sigma_{\rm
ess}(\Delta^p_M)=
\emptyset,$$ whilst if $p\in \left\{0,1,n-1,n\right\}$ $$\sigma_{\rm
ess}(\Delta^p_M)= [0,+\infty).$$ \end{enumerate}
\end{thm}


\end{document}